\documentclass[twoside,12pt]{article}
\usepackage{amsmath}
\usepackage{amssymb}
\usepackage{mathrsfs}

\usepackage{cite}

\def\sign{\mathop{\rm sign}}

\def\Sign{\mathop{\rm Sign}}

\title{Solvability of a class of phase field systems \\ related to a sliding mode control problem}

\author{Michele Colturato\\
Dipartimento di Matematica, Universit\`a degli Studi di Pavia\\
Via Ferrata~1, 27100 Pavia, Italy\\
E-mail: \texttt{michele.colturato01@universitadipavia.it}}
 
\date{}

\newcommand\testopari{\sc Michele Colturato}
\newcommand\testodispari{\sc Solvability of a phase field system related to a SMC problem}
\newtheorem{Rem1}{Remark}[section]
\newtheorem{Rem2}[Rem1]{Remark}
\newtheorem{Teo-esistenza}[Rem1]{Theorem - (Existence)}
\newtheorem{Teo3}[Rem1]{Theorem - (Uniqueness and continuous dependence)}
\newtheorem{Rem3}[Rem1]{Remark}
\newtheorem{Rem11}[Rem1]{Remark}
\newtheorem{Rem4}[Rem1]{Remark}

\begin{document}
\maketitle

\begin{abstract}
We consider a phase-field system of Caginalp type perturbed by the presence of an additional maximal monotone nonlinearity.
Such a system arises from a recent study of a sliding mode control problem.
We prove existence of strong solutions. Moreover, under further assumptions, we show the continuous dependence on the initial data
and the uniqueness of the solution.

\vspace{2mm}

\vspace{2mm}
\noindent \textbf{AMS (MOS) subject clas\-si\-fi\-ca\-tion:} 35K61, 35K25, 
35B25, 35D30, 80A22.

\end{abstract}

\section{Introduction}
In the present contribution we consider the phase-field system
\begin{equation} \label{originale1111}
\partial_t(\mu  + \ell \phi) -k \Delta \mu  + \zeta = f \ \ \textrm{a.e. in $Q:= (0,T) \times \Omega$,} 
\end{equation}
\begin{equation} \label{originale21111}
\partial_t \phi - \nu \Delta \phi + \xi + \pi(\phi) = \gamma \mu   \ \ \textrm{a.e. in $Q$,}
\end{equation}
\begin{equation} 
\zeta(t) \in A(\mu (t) + \alpha \phi(t) - \Xi^*) \ \textrm{for a.e. $t \in (0,T)$,}  
\end{equation}
\begin{equation} 
\xi \in \chi(\phi) \ \textrm{a.e. in $Q$,}
\end{equation}
whe  re   $\Omega$ is the   domain in which the   e  volution take  s place  ,
$T$~is some   final time  ,
$\mu $~de  note  s the   re  lative   te  mpe  rature   around some   critical value  
that is take  n to be   $0$ without loss of ge  ne  rality,
and $\phi$ is the   orde  r parame  te  r.
More  ove  r, $\ell$, $k$, $\nu$,  $\gamma$ and $\alpha$ are   positive   constants,  $\Xi^*$ is a function
in $H^2(\Omega)$ with null outward normal de  rivative   on the   boundary of $\Omega$ and
$f$~is a source te  rm. The above   syste  m is comple  me  nte  d by
homoge  ne  ous Ne  umann boundary conditions for both $\mu $ and~$\phi$, that~is,
\begin{equation}
\partial_{N}\mu  = 0, \quad \quad \partial_{N}\phi = 0 \quad \textrm{on $\Sigma:= (0,T) \times \Gamma$},
\end{equation}
whe  re   $\Gamma$ is the   boundary of $\Omega$ and $\partial_{N}$ is the   outward normal de  rivative  ,
and by the   initial conditions 
\begin{equation} \label{originale51111}
\mu (0)=\mu _0, \quad \quad \phi(0)=\phi_0  \quad  \textrm{in $\Omega$.}
\end{equation} 
The   te  rm $\xi + \pi(\phi)$, appe  aring in \eqref{originale21111}, re  pre  se  nts the   de  rivative   of a double  -we  ll pote  ntial $F$
de  fine  d as the sum 
\begin{equation} 
F= \widehat{\chi} + \widehat{\pi},
\end{equation}
whe  re
	\begin{equation} 
\widehat{\chi}: \mathbb{R} \longrightarrow [0, + \infty] \textrm{ is proper, l.s.c. and convex with $\widehat{\chi}(0)=0$,} 
		\end{equation}
		\begin{equation} 
\widehat{\pi}:  \mathbb{R} \rightarrow \mathbb{R}, \ \textrm{$\widehat{\pi} \in C^1(\mathbb{R})$ with $ \pi := \widehat{\pi}'$ Lipschitz continuous.}
	\end{equation}
Since   $\widehat{\chi}$ is prope  r, l.s.c. and conve  x, the   subdiffe  re  ntial $\partial \widehat{\chi} =: \chi$ is well de  fine  d and is a maximal monotone   graph.
For a compre  hensive discussion of the the  ory of 
maximal monotone   ope  rators, we refe  r, e.g., to \cite{Barbu, Brezis, Show}.	
In our proble  m we   also conside  r a maximal monotone   ope  rator
\begin{equation} \label{A111111}
A: \ H:=L^2(\Omega) \longrightarrow  2^H
\end{equation}
such that $0 \in A(0)$ and 
\begin{equation}  \label{A1111112}
\| v \|_H \leq C (1 + \| x \|_H )  \quad \quad \textrm{  for all} \ x \in H, \ v \in A x,
\end{equation}
for some constant $C>0$.  

The   proble  m \eqref{originale1111}--\eqref{originale51111} under study is an interesting de  ve  lopme  nt of the   following
simple   ve  rsion of the   phase  -fie  ld syste  m of Caginalp type   (se  e~\cite{Cag}):
\begin{equation} \label{Iprima}
\partial_t(\mu  + \ell \phi) -k \Delta \mu  = f  \quad  \textrm{in $Q$,} 
\end{equation}
\begin{equation} \label{Iseconda} 
\partial_t \phi - \nu \Delta \phi + F'(\phi) = \gamma \mu   \quad  \textrm{in $Q$.} 
\end{equation}
As already noticed, $F' \cong \xi + \pi$ is related to a double-well potential~$F$.
Typical examples for $F$ are
\begin{equation} \label{regpot}
F_{reg}(r) = \frac{1}{4}(r^2-1)^2, \quad \quad r \in \mathbb{R},
\end{equation}
\begin{equation} \label{logpot}
F_{log}(r) = ((1+r) \ln{(1+r)} + (1-r) \ln{(1-r)})-c_0 r^2, \quad r \in (-1,1),
\end{equation}
\begin{equation} \label{obspot}
F_{obs}(r)= I(r) -c_0 r^2, \quad \quad r \in \mathbb{R},
\end{equation} 
whe  re   $c_0>1$ in \eqref{logpot} in order to produce a double   we  ll,
while   $c_0$ is an arbitrary positive   numbe  r in~\eqref{obspot},
and the function $I$ in~\eqref{obspot} is the indicator function of~$[-1,1]$, i.e.,
it take  s the   value  s $0$ or $+\infty$ according to whe  the  r or not $r$ be  longs to~$[-1,1]$.
The   pote  ntials \eqref{regpot} and \eqref{logpot} are 
the   usual classical re  gular pote  ntial and the   so-calle  d logarithmic pote  ntial, re  spe  ctively.

The  we ll-pose dne ss, the  long-time  be havior of solutions, 
and also the  re late d optimal control proble ms 
conce rning Caginalp-type  syste ms have  be e n wide ly studie d in the  lite rature . 
We  re fe r, without any sake  of comple te ne ss, e .g., to \cite{BrokSpr, EllZheng, GraPetSch, KenmNiez, Lau} 
and refe re nce s the re in for the  we ll-pose dne ss and long time  be havior re sults and to 
\cite{CGM, CoGiMaRo, CoGiMaRorrrrrrrrr, HoffJiang, HKKY} for 
the  tre atme nt of optimal control proble ms.

The  pape r \cite{BaCoGiMaRo} is re late d to control proble ms, but it goe s in the  dire ction
of de signing sliding mode  controls (SMC) for a particular phase -fie ld syste m. The  main
obje ctive  of the  authors is to find some  state -fe e dback control laws $(\mu , \phi) \mapsto u(\mu , \phi)$ that can be 
that, once  inse rte d
into the  e quations, force  the  solution to re ach some  submanifold of the  phase  space , in
finite  time , the n slide  along it.
The  first analytical difficulty consists in de riving the  e quations gove rning the  sliding
mode s and the  conditions for this motion to e xist. The  proble m ne e ds the  de ve lopme nt
of spe cial me thods, since  the  conve ntional the ore ms re garding e xiste nce  and unique ne ss
of solutions are  not dire ctly applicable . More ove r, the  authors ne e d to manipulate  the  syste m
through the  control in orde r to constrain the  e volution on the  de sire d sliding manifold.

In particular, in the  pape r \cite{BaCoGiMaRo} the authors conside r the  ope rator $\Sign : H \longrightarrow 2^H$ define d as
$\Sign(v) = \frac{v}{\| v \|}$, if $v \neq 0$ and $\Sign(0) = B_1(0)$, if $v = 0$, where $B_1(0)$ is the  close d unit ball of $H$.
$\Sign$ is a maximal monotone  ope rator on $H$ 
and is a nonlocal counte rpart of the  ope rator  
$\sign : \mathbb{R} \longrightarrow 2^{\mathbb{R}}$ de fine d as
$\sign (r) = \frac{r}{| r |}$, if $r\neq 0$  and $\sign (0) = [ - 1, 1 ]$, if $r = 0$.
The n the  authors of \cite{BaCoGiMaRo} de al with the  syste m
\begin{equation} \label{barbucolligilardi1}
\partial_t(\mu  + \ell\phi) -k \Delta \mu  = f- \rho \sigma  \ \ \textrm{a.e. in $Q$,} 
\end{equation}
\begin{equation} \label{barbucolligilardi2}
\partial_t \phi - \nu \Delta \phi + \xi + \pi(\phi) = \gamma \mu   \ \ \textrm{a.e. in $Q$,}
\end{equation}
\begin{equation} \label{barbucolligilardi4}
\sigma(t) \in \Sign(\mu (t) + \alpha \phi(t) - \Xi^*) \ \ \textrm{for a.e. $t \in (0,T)$},
\end{equation}
\begin{equation} 
\xi \in \chi(\phi) \ \ \textrm{a.e. in $Q$,}
\end{equation}
\begin{equation}
\partial_{N}\mu  = 0, \quad \quad \partial_{N}\phi = 0 \quad \textrm{on $\Sigma$},
\end{equation}
\begin{equation} \label{barbucolligilardi3}
\mu (0)=\mu _0, \quad \quad \phi(0)=\phi_0 \quad  \textrm{in $\Omega$,} 
\end{equation} 
which turns out to be  a particular case  of \eqref{originale1111}--\eqref{originale51111} with $A = \rho \Sign$.
The  pape r \cite{BaCoGiMaRo} is mostly conce rne d with the  sliding mode  prope rty for \eqref{barbucolligilardi1}--\eqref{barbucolligilardi3}.
In this contribution we  de al with \eqref{originale1111}--\eqref{originale51111}, which turns out to be a particular
ge ne ralization of the  proble m \eqref{barbucolligilardi1}--\eqref{barbucolligilardi3}
since  we  only re quire  \eqref{A111111}--\eqref{A1111112} for the  maximal monotone  ope rator $A$.
We  prove  e xiste nce  and re gularity of the  solutions for the  problem \eqref{originale1111}--\eqref{originale51111},
as well as the  unique ne ss  and the  continuous de pe nde nce  on the  initial data in case  $\alpha = \ell$.
In orde r to obtain our re sults, we  first make  a change  of variable . We  se t:
\begin{equation} 
\Xi = \mu  + \alpha \phi - \Xi^*.
\end{equation}
Consequently, the previous system \eqref{originale1111}--\eqref{originale51111} becomes
\begin{equation} \label{introduzione1}
\partial_t(\Xi + (\ell-\alpha) \phi) - k \Delta \Xi + k \alpha \Delta \phi + \zeta= f- k \Delta \Xi^{*} \ \ \textrm{a.e. in $Q$,} 
\end{equation}
\begin{equation} 
\partial_t \phi - \nu \Delta \phi + \xi + \pi(\phi) = \gamma (\Xi - \alpha \phi + \Xi^{*})  \ \ \textrm{a.e. in $Q$,}
\end{equation}
\begin{equation} 
\zeta(t) \in A(\Xi(t)) \ \ \textrm{for a.e. $t \in (0,T)$},
\end{equation}
\begin{equation} \label{introduzione5}
\xi \in \chi(\phi) \ \ \textrm{a.e. in $Q$.}
\end{equation}
\begin{equation} 
\partial_{N}\Xi = 0, \quad \quad \partial_{N}\phi = 0  \ \ \textrm{on $\Sigma$},
\end{equation} 
\begin{equation} \label{introduzione99}
\Xi(0) = \Xi_0, \quad \quad \phi(0) = \phi_0 \quad \textrm{in $\Omega$.}
\end{equation} 
From now on, 
we  re fe r to the  initial
and boundary value  proble m \eqref{introduzione1}--\eqref{introduzione99} as Problem $(P)$.
In orde r to prove  the  e xiste nce  of solutions, we  first conside r the  approximating proble m $(P_{\varepsilon})$,
obtaine d from proble m $(P)$ by approximating $A$ and $\chi$ by the ir Yosida
re gularizations. The n we  construct a furthe r approximating proble m  $(P_{\varepsilon, n})$, obtaine d from $(P_{\varepsilon})$
by a Faedo-Galerkin sche me  base d on a syste m of e ige nfunctions $\{ v_n \} \subseteq W $, whe re
\begin{equation}
W = \{ u \in H^2(\Omega): \ \partial_{N} u = 0 \ \textrm{on} \ \partial \Omega \}.
\end{equation}
The n, we  prove  the  e xiste nce  of a local solution for $(P_{\varepsilon, n})$ and,  passing to the  limit as $n \rightarrow + \infty$,
we  infe r that the  limit of some  subse que nce  of solutions for $(P_{\varepsilon, n})$ yie lds a solution of $(P_{\varepsilon})$.
Finally, we  pass to the  limit as $\varepsilon \searrow 0$
and show that some limit of a subse que nce  yie lds a solution of $(P)$.

Ne xt, we  le t $\alpha = \ell$ and  write  proble m $(P)$ 
for two diffe re nt se ts of initial data $f_i$, $\Xi^*_i$, $\Xi_{0_{i}}$ and $\phi_{0_{i}}$, $i=1,2$.
By  pe rforming suitable  contracting e stimate s for the  diffe re nce  of the  corre sponding solutions,
we  de duce  the  continuous de pe nde nce  re sult 
whe nce  the  unique ne ss prope rty is also achie ve d.

\section{Main results}
\setcounter{equation}{0}

\subsection{Preliminary assumptions}
We  assume  $\Omega \subseteq \mathbb{R}^3$ to be  a bounde d domain of class $C^1$ and we  write  $| \Omega |$ for its Le be sgue  me asure.
More over, $\Gamma$ and $\partial_{N}$ still stand for the  boundary of $\Omega$ and the  outward normal de rivative , re spe ctive ly.
Give n a finite  final time  $T > 0$, for e very $t \in (0,T]$ we set
\begin{equation}
Q_t = (0,t) \times \Omega, \ \ Q = Q_T,  
\end{equation}
\begin{equation}
\Sigma_t = (0,t) \times \Gamma, \ \ \Sigma = \Sigma_T. 
\end{equation}
In the following, we set for brevity:
\begin{equation}
H = L^2(\Omega), \ \ \ \ V = H^1(\Omega), \ \ \ \ V_0 = H^1_0(\Omega),
\end{equation}
\begin{equation} \label{W}
W = \{ u \in H^2(\Omega): \ \partial_{N} u = 0 \ \textrm{on} \ \partial \Omega \},
\end{equation}
with usual norms $\| \cdot \|_{H}$, $\|\cdot \|_{V}$
and inner products $(\cdot,\cdot )_{H}$, $(\cdot ,\cdot )_{V}$, respectively. 
Now we  de scribe  the  proble m unde r conside ration.
We assume that 
\begin{equation} \label{parametri}
\ell, \ \alpha, \ k, \ \nu, \ \gamma \in (0, + \infty),
\end{equation}
\begin{equation} \label{f}
f \in L^2(Q),
\end{equation}
\begin{equation} \label{eta}
\Xi^{*} \in W, 
\end{equation}
\begin{equation} \label{0}
\Xi_0, \ \phi_0 \in V,
\end{equation}
\begin{equation} \label{betal10}
\widehat{\chi}(\phi_0) \in L^1(\Omega).
\end{equation}
We  introduce  the  double -we ll pote ntial $F$ as the sum 
\begin{equation} 
F= \widehat{\chi} + \widehat{\pi},
\end{equation}
where
	\begin{equation} \label{beta}
\widehat{\chi}: \mathbb{R} \longrightarrow [0, + \infty] \textrm{ is proper, l.s.c. and convex with $\widehat{\chi}(0)=0$,} 
		\end{equation}
\begin{equation} \label{pi}
\widehat{\pi}:  \mathbb{R} \rightarrow \mathbb{R}, \ \textrm{$\widehat{\pi} \in C^1(\mathbb{R})$ with $ \pi := \widehat{\pi}'$ Lipschitz continuous.}
	\end{equation}
Since $\widehat{\chi}$ is proper, lower semicontinuous and convex, the subdifferential $\partial \widehat{\chi} =: \chi$ is well defined. We denote by 
$D(\chi)$ and $D(\widehat{\chi})$ the effective domains of $\chi$ and $\widehat{\chi}$, respectively. Thanks to these assumptions, $\chi$ is a maximal monotone graph.	Moreover, as $\widehat{\chi}$ takes on its minimum in $0$, we have that 
$0 \in \chi(0)$.

\begin{Rem1}
\emph{We introduce the operator $\mathcal{B}$ induced by $\chi$ on $L^2(Q)$ in the following way:}
\begin{equation} \label{betagrande1}
\mathcal{B}: L^2(Q) \longrightarrow L^2(Q)
\end{equation}
\begin{equation} \label{betagrande2}
\xi  \in \mathcal{B}(\phi) \Longleftrightarrow \xi(x,t) \in \chi(\phi(x,t)) \quad \textrm{\emph{ for a.e.} $(x,t) \in Q$.}
\end{equation}
\emph{We notice that}
\begin{equation} 
\chi= \partial \widehat{\chi},  \quad \quad \quad \quad \mathcal{B} = \partial \Phi, 
\end{equation}
\emph{where} 
\begin{equation} 
\Phi: \  L^2(Q) \longrightarrow (-\infty , + \infty]
\end{equation}
\begin{equation}
\Phi(u)=
\left\{ \begin{array}{ll}
\int_Q{\widehat{\chi}(u)} & \textrm{\emph{if $u \in L^2(Q)$ and $\widehat{\chi}(u) \in L^1(Q)$},} \\
+ \infty & \textrm{\emph{elsewhere, with  $u \in L^2(Q)$.}}
\end{array}
\right.
\end{equation}	
\end{Rem1}
		
\paragraph{The maximal monotone operator $A$.}
In our problem a maximal monotone operator
\begin{equation} \label{A1}
A: H \longrightarrow  H
\end{equation}
also appears. We assume that
\begin{equation} \label{A2}
0 \in A(0)
\end{equation}
and that there exists a constant $C>0$ such that
\begin{equation} \label{stimaA}
\| v \|_H \leq C (1 + \| \Xi \|_H )  \quad \quad \hbox{  for all} \ \Xi \in H, \ v \in A \Xi.
\end{equation}

\begin{Rem2}
\emph{We introduce the operator $\mathcal{A}$ induced by $A$ on $L^2(0,T;H)$ in the following way}
\begin{equation} \label{Agrande1}
\mathcal{A}: L^2(0,T;H) \longrightarrow L^2(0,T;H)
\end{equation}
\begin{equation} \label{Agrande2}
\zeta  \in \mathcal{A}(\Xi) \Longleftrightarrow \zeta(t) \in A(\Xi(t)) \quad \textrm{\emph{ for a.e.} $t \in (0,T)$.}
\end{equation}
\emph{We  notice that $\mathcal{A}$ is a maximal monotone operator.}
\end{Rem2}

\subsection{Examples of operators $A$}
Now, we provide some examples
of maximal monotone operators fulfilling our assumptions.
\paragraph{Example 1.}
We consider the operator
\begin{equation} \label{segno111111111111144444444444444444444444}
\sign \ : \mathbb{R} \longrightarrow 2^{\mathbb{R}}
\end{equation}
\begin{equation} \label{segno1111111111111}
\sign(r) =  
\left\{ \begin{array}{ll}
\frac{r}{| r |}  								& \textrm{if $r \neq 0$}, \\
\textrm{$[ - 1, 1 ]$}  										& \textrm{if $ r =0  $.}
\end{array}
\right.
\end{equation}
Notice that $\sign$ induces a maximal monotone operator on $H$.
\paragraph{Example 2.}
We define  the operator $\Sign$ as the nonlocal counterpart of the operator $\sign$ (see \eqref{segno111111111111144444444444444444444444}--\eqref{segno1111111111111}):
\begin{equation} \label{Sign}
\Sign: \quad \quad H \longrightarrow 2^H
\end{equation}
\begin{equation}
\Sign(v) = 
\left\{ \begin{array}{ll}
\frac{v}{\| v \|} 								& \textrm{if $v \neq 0$}, \\
B_1(0)													& \textrm{if $ v =0  $,}
\end{array}
\right.
\end{equation}
where $B_1(0)$ is the closed unit ball of $H$.
$\Sign$ is the subdifferential of the map $ \| \cdot \| : H \rightarrow  \mathbb{R} $ and
is a maximal monotone operator on $H$ which satisfies \eqref{A2}--\eqref{stimaA}.
\paragraph{Example 3.}
We consider the operator
	\begin{equation}
A_1: \ \mathbb{R} \longrightarrow \mathbb{R}
\end{equation}
\begin{equation}
A_1(r) = 
\left\{ \begin{array}{ll}
\alpha_1 r  								& \textrm{if $r<0$}, \\
0  													& \textrm{if $0 \leq r \leq 1 $}, \\
\alpha_2 r  		& \textrm{if $r > 1 $}, \\
\end{array}
\right.
\end{equation}
where $\alpha_1$ and $\alpha_2 $ are positive coefficients.
We observe that  $A_1$ is a maximal monotone operator on $\mathbb{R}$,
whose graph consists of an horizontal line segment and
two rays of slope $\alpha_1$, $\alpha_2 $. Moreover, $0 \in A_1(0)$ and
\begin{equation}
| v | \leq C (1 + | r | )  \quad \quad \hbox{  for all} \ r \in \mathbb{R}, \ v \in A_1(r),
\end{equation}
with $C = \max{(\alpha_1, \alpha_2 )}$. Then $A_1$ satisfies \eqref{A2}--\eqref{stimaA}.
We notice that $A_1$ corresponds to the graph which correlates the enthalpy to the temperature in the Stefan problem (see, e.g., \cite{Dam77, duvaut, Fri68}).
\paragraph{Example 4.}
We consider the operator
\begin{equation}
A_2: \ H \longrightarrow H
\end{equation}
\begin{equation}
A_2(v) = \alpha |v|^{q-1} v, 
\end{equation}
where  $0<q<1$ and $\alpha $ is a function in $ L^{\infty}(\Omega)$ with $\alpha(x) \geq 0$ for a.e. $x \in \Omega$.
We observe that $A_2$ induces a (nonlocal) multivalued  maximal monotone operator on $H$, with $0 \in A_2(0)$.
Moreover, $A_2$ can be considered a weighted perturbation of the operator appearing
in the porous media equation and in the fast diffusion equation (see, e.g.,  \cite{Dib83, Hui07, Vaz07}).


\subsection{Setting of the problem and results}
Now, we state the problem under consideration. 
We look for a pair $(\Xi,\phi)$ satisfying at least the regularity requirements 
\begin{equation} \label{regloarita}
\Xi, \ \phi \ \in H^1(0,T; H) \cap L^{\infty}(0,T; V) \cap L^2(0,T; W) 
\end{equation}
and solving the problem $(P)$:
\begin{equation} \label{iniziale1}
\partial_t(\Xi + (\ell-\alpha) \phi) -k \Delta \Xi + k \alpha \Delta \phi + \zeta = f-k \Delta \Xi^{*} \ \ \textrm{a.e. in $Q$,} 
\end{equation}
\begin{equation} \label{iniziale2}
\partial_t \phi - \nu \Delta \phi + \xi + \pi(\phi) = \gamma (\Xi - \alpha \phi + \Xi^{*})  \ \ \textrm{a.e. in $Q$,}
\end{equation}
\begin{equation} \label{iniziale3}
\zeta(t) \in A(\Xi(t)) \ \ \textrm{for a.e. $t \in (0,T)$,}  
\end{equation}
\begin{equation} \label{iniziale4}
\xi \in \chi(\phi) \ \ \textrm{a.e. in $Q$,}
\end{equation}
\begin{equation} \label{iniziale4-aggiunta}
\partial_{N}\mu  = 0, \quad \quad \partial_{N}\phi = 0 \quad \textrm{on $\Sigma$},
\end{equation}
\begin{equation} \label{iniziale5}
\Xi(0) = \Xi_0, \quad \quad \phi(0) = \phi_0 \quad \textrm{in $\Omega$.}
\end{equation}
We notice that the homogeneous Neumann boundary conditions for both $\Xi$ and $\phi$ required by \eqref{iniziale4-aggiunta}
follow from \eqref{regloarita}, due to the definition of $W$ (see \eqref{W}).

\begin{Teo-esistenza} \label{Teo-esistenza}
Assume \eqref{parametri}--\eqref{betal10}, \eqref{beta}--\eqref{pi} and \eqref{A1}--\eqref{stimaA}. 
Then problem $(P)$ (see \eqref{iniziale1}--\eqref{iniziale5})
has at least a solution $(\Xi,\phi)$ satisfying the regularity requirements \eqref{regloarita}. 
\end{Teo-esistenza}

\begin{Teo3} \label{Teo3}
Assume \eqref{parametri}--\eqref{betal10}, \\ \eqref{beta}--\eqref{pi} and \eqref{A1}--\eqref{stimaA}. 
If $\alpha = \ell$, the solution $(\phi, \Xi)$ of problem $(P)$ (see \eqref{iniziale1}--\eqref{iniziale5}) is unique. 
Moreover, if $f_i$, $\Xi^{*}_i$, $\Xi_{0_i}$, $\phi_{0_i}$, $i=1,2$, 
are given as in \eqref{f}--\eqref{0}
and $(\phi_i, \Xi_i)$, $i=1,2$,  are the corresponding solutions,
then the estimate
\begin{equation} \nonumber
\| \Xi_1 - \Xi_2 \|_{L^{\infty}(0,T; H) \cap L^{2}(0,T; V) } + \|  \phi_1 - \phi_2 \|_{L^{\infty}(0,T; H) \cap L^{2}(0,T; V)} 
\end{equation}
\begin{equation}
\leq C ( \|f_1 -f_2 \|_{L^{2}(Q)} + \| \Xi^{*}_1 - \Xi^{*}_2 \|_W + \| \Xi_{0_1} - \Xi_{0_2} \|_H + \| \phi_{0_1} - \phi_{0_2} \|_H )
\end{equation}
holds true for some constant $C$ that depends only on $\Omega$, $T$
and the parameters $\ell$, $\alpha$, $k$, $\nu$, $\gamma$.
\end{Teo3}

\section{Proof of the existence theorem}
\setcounter{equation}{0}
This section is devoted to the proof of Theorem \ref{Teo-esistenza}.

\subsection{The approximating problem $(P_{\varepsilon})$} 
\paragraph{Yosida regularization of $A$.} 
We introduce the Yosida regularization of $A$. For $\varepsilon>0$ we define
\begin{equation} \label{richiamoA1}
A_{\varepsilon}: H \longrightarrow  H, \quad \quad
A_{\varepsilon} = \frac{I - (I + \varepsilon A)^{-1}}{\varepsilon},
\end{equation} 
where $I$ denotes the identity operator.
Note that $A_{\varepsilon}$ is is Lipschitz-continuous (with Lipschitz
constant $\frac{1}{\varepsilon}$), maximal monotone, 
and satisfies the following properties. Denoting by $J_{\varepsilon}= (I + \varepsilon A)^{-1}$ the resolvent operator, for all $\delta > 0$ we have that
	\begin{equation} \label{AJe}
	A_{\varepsilon}\Xi \in A(J_{\varepsilon}\Xi),
	\end{equation}
	\begin{equation}
	(A_{\varepsilon})_\delta = A_{\varepsilon + \delta},
	\end{equation}
	\begin{equation}
	\| A_{\varepsilon}\Xi \|_H \leq \| A^0 \Xi \|_H,
	\end{equation}
	\begin{equation}
	\lim_{\varepsilon \rightarrow 0} \| A_{\varepsilon}\Xi \|_H = \| A^0 \Xi \|_H,
	\end{equation}
where $A^0 \Xi$ is the element	of the image of $A$ having minimal norm. 

\begin{Rem3}
\emph{We point out a key property of} $A_{\varepsilon}$\emph{, 
which is a consequence of \eqref{stimaA}. There exists a positive constant $C$, 
independent of $\varepsilon$, such that}
\begin{equation} \label{Ae}
\| A_{\varepsilon}\Xi \|_H \leq C (1 + \| \Xi \|_H ) \quad \quad \hbox{\emph{ for all}} \ \Xi \in H, \ v \in A \Xi.
\end{equation}
\emph{Indeed, notice that} $0 \in A(0)$ \emph{and} $0 \in I(0)$\emph{: consequently, for every} $\varepsilon >0$\emph{,} $0 \in (I + \varepsilon A)(0)$. \emph{This fact implies that }
$ J_\varepsilon (0)=0$.
\emph{Moreover, since } $A$ \emph{is  maximal monotone, } $J_\varepsilon$ \emph{is a contraction. Then, from \eqref{stimaA} and \eqref{AJe}, it follows that}
\begin{equation} \nonumber
\begin{array}{ll}
\| A_{\varepsilon}\Xi \|_H & \leq C (\| J_{\varepsilon}\Xi\|_H  + 1)\\
& \leq C (\| J_{\varepsilon}\Xi - J_{\varepsilon}0 \|_H + \| J_{\varepsilon}0 \|_H + 1)\\
& \leq C (\| \Xi \|_H + 1).
\end{array}													
\end{equation}
\end{Rem3}

\paragraph{Yosida regularization of $\chi$.} 
We introduce the Yosida regularization of $\chi$. For $\varepsilon>0$ we define
\begin{equation} \label{richiamobeta1}
\chi_{\varepsilon}: \mathbb{R} \longrightarrow \mathbb{R},
\quad \quad
\chi_{\varepsilon} = \frac{I - (I + \varepsilon \chi)^{-1}}{\varepsilon}. 
\end{equation}
We remark that $\chi_{\varepsilon}$ is Lipschitz continuous (with Lipschitz constant $\frac{1}{\varepsilon}$) and satisfies the following properties. Denoting by $R_{\varepsilon}= (I + \varepsilon \chi)^{-1}$ the resolvent operator, for all $\delta > 0$ and for every $\phi \in D(\chi)$ we have that
	\begin{equation}
	\chi_{\varepsilon}(\phi) \in \chi(R_{\varepsilon}\phi),
  \end{equation}
		\begin{equation}
	(\chi_{\varepsilon})_\delta = \chi_{\varepsilon + \delta}, 
	\end{equation}
	\begin{equation}
	|\chi_{\varepsilon}(\phi)| \leq |\chi^0(\phi)|,
	\end{equation}
	\begin{equation}
 	\lim_{\varepsilon \rightarrow 0} \chi_{\varepsilon}(\phi) =  \chi^0(\phi), 
	\end{equation}
where $\chi^0(\phi)$ is the element	of the image of $\chi$ having minimal modulus.

\paragraph{Regularization of $\widehat{\chi}$.} 
We introduce the Moreau-Yosida regularization of $\widehat{\chi}$. For $\varepsilon>0$ we define
\begin{equation}
\widehat{\chi}_{\varepsilon}: \mathbb{R} \longrightarrow [0, +\infty],
\quad \quad \widehat{\chi}_{\varepsilon} = \frac{I - (I + \varepsilon \widehat{\chi})^{-1}}{\varepsilon}. 
\end{equation}
We recall that
	\begin{equation} \label{questaservepoi}
	\widehat{\chi}_{\varepsilon}(\phi) \leq \widehat{\chi}(\phi) \quad \textrm{for every $\phi  \in D(\widehat{\chi})$.}
	\end{equation}
We also observe that $\chi_{\varepsilon}$ is the Fr\emph{$\acute{e}$}chet derivative of $\widehat{\chi}_{\varepsilon}$. Then, for every $\phi_1, \phi_2  \in D(\widehat{\chi})$, we have that
\begin{equation} \label{questaservepoimolto}
\widehat{\chi}_{\varepsilon}(\phi_2) = \widehat{\chi}_{\varepsilon}(\phi_1) + \int_{\phi_1}^{\phi_2} \chi_{\varepsilon}(s) \ ds.
\end{equation}

\paragraph{Approximating problem $(P_{\varepsilon})$.}
We denote by $f_{\varepsilon}$ a regularization of $f$ constructed in such a way that
\begin{equation} \label{convergenzadif}
f_{\varepsilon} \in C^1([0,T]; H) \textrm{ for all $\varepsilon >0$}, \quad \quad
\lim_{\varepsilon \rightarrow 0} \| f_{\varepsilon} - f \|_{L^2(0,T; H)} = 0.
	\end{equation}
Then, we look for a pair $(\Xi_{\varepsilon},\phi_{\varepsilon})$ satisfying at least the regularity requirements 
\begin{equation} \label{regloarita2}
\Xi_{\varepsilon}, \ \phi_{\varepsilon} \ \in H^1(0,T; H) \cap L^{\infty}(0,T; V) \cap L^2(0,T; W)
\end{equation}
and solving the approximating  problem $(P_{\varepsilon})$:
\begin{equation} \label{equazionedelsecondolimite1}
\partial_t(\Xi_{\varepsilon} + (\ell-\alpha) \phi_{\varepsilon}) -k \Delta \Xi_{\varepsilon} + k \alpha \Delta \phi_{\varepsilon} + \zeta_{\varepsilon} = f_{\varepsilon}-k \Delta \Xi^{*} \ \ \textrm{a.e. in $Q$,}
\end{equation}
\begin{equation}
\partial_t \phi_{\varepsilon} - \nu \Delta \phi_{\varepsilon} + \xi_{\varepsilon} + \pi(\phi_{\varepsilon}) = \gamma (\Xi_{\varepsilon} - \alpha \phi_{\varepsilon} + \Xi^{*}) \ \ \textrm{a.e. in $Q$,}
\end{equation}
\begin{equation}
\zeta_{\varepsilon}(t) = A_{\varepsilon} \Xi_{\varepsilon}(t) \ \ \textrm{for a.e. $t \in (0,T)$,}
\end{equation}
\begin{equation}
\xi_{\varepsilon} = \chi_{\varepsilon}(\phi_{\varepsilon}) \ \ \textrm{a.e. in $Q$,}
\end{equation}
\begin{equation} \label{utilinaqui1}
\partial_{N}\Xi_{\varepsilon} = 0, \quad \quad \partial_{N}\phi_{\varepsilon} = 0 \quad \textrm{on $\Sigma$},
\end{equation}
\begin{equation} \label{equazionedelsecondolimite2}
\Xi_{\varepsilon}(0) = \Xi_{0},  \quad \quad  \phi_{\varepsilon}(0) = \phi_{0} \quad  \textrm{in $\Omega$,}
\end{equation}
where $A_{\varepsilon}$ and $\chi_{\varepsilon}$  are the Yosida regularizations of   $A$ and $\chi$
defined in   \eqref{richiamoA1} and \eqref{richiamobeta1}, respectively.
We notice that the homogeneous Neumann boundary conditions for both $\Xi_{\varepsilon}$ and $\phi_{\varepsilon}$ 
required by \eqref{utilinaqui1}
follow from \eqref{regloarita2} due to the definition of $W$ (see \eqref{W}).

\begin{Rem11}
\emph{We can define} $f_{\varepsilon}$ \emph{as the regularization of} $f$ \emph{obtained solving}
\begin{equation}
\left\{ \begin{array}{ll}
-\varepsilon f''_{\varepsilon}(t) + f_{\varepsilon}(t)  = f(t), \quad t \in (0,T),  \\
f_{\varepsilon}(0) = f_{\varepsilon}(T) = 0. \\
\end{array}
\right.
\end{equation}
\emph{Thanks to Sobolev immersions and elliptic regularity, we obtain that}
\begin{equation}
f_{\varepsilon} \in C^1([0,T]; H) \quad \textrm{\emph{and}} \quad
\lim_{\varepsilon \rightarrow 0} \| f_{\varepsilon} - f \|_{L^2(0,T; H)} = 0.
\end{equation}
\end{Rem11}


\subsection{The approximating problem $(P_{\varepsilon,n})$}
Now, we apply the Faedo-Galerkin method to the approximating problem $(P_{\varepsilon})$.
We consider the orthonormal basis $\{ v_i \}_{\{ i \geq 1 \}}$
of $V$ formed by the normalized eigenfunctions of the Laplace operator with homogeneous Neumann boundary condition, that is 
\begin{equation} \label{baseortonormale}
\left\{ \begin{array}{ll}
- \Delta v_i = \lambda_i v_i & \textrm{in $\Omega$,} \\
\partial_{N} v_i = 0 & \textrm{on $\partial \Omega$.} \\
\end{array}
\right.
\end{equation}
Note that, owing to the regularity of $\Omega$, $ v_i \in W$ for all $i \geq 1$.
Then, for any integer $n \geq 1$, we denote by $V_n$ the $n$-dimentional subspace of $V$ spanned by $\{ v_1 , \cdots , v_n \}$.
Hence, $ \{ V_n \} $ is a sequence of finite dimensional subspaces such that $ \bigcup_{n=1}^{+\infty} V_n$ is dense in $V$
and $V_k \subseteq V_n$ for all $k \leq n$.

\paragraph{Definition of the approximating problem $(P_{\varepsilon,n})$.}
We first approximate the initial data $\Xi_0$ and $\phi_0$. We set
\begin{equation} 
\Xi_{0,n} = P_{V_n} \Xi_0, \quad \quad 
\phi_{0,n} = P_{V_n} \phi_0. 
\end{equation}
We notice that
\begin{equation} \label{datiinizialifaedo}
\lim_{n \rightarrow + \infty} \| \Xi_{0,n} - \Xi_0 \|_V = 0 \quad \textrm{and} \quad
\lim_{n \rightarrow + \infty} \| \phi_{0,n} - \phi_0 \|_V = 0. 
\end{equation}
Note that the convergence provided by \eqref{datiinizialifaedo} assures that $\Xi_{0,n}$ and $\phi_{0,n}$ are bounded in $V$.
Now, we introduce the new approximating problem $(P_{\varepsilon,n})$. We look for $t_n \in ]0, T]$
and a pair $(\Xi_{\varepsilon,n},\phi_{\varepsilon,n})$ 
(in the following we will write $(\Xi_n,\phi_n)$ instead of $(\Xi_{\varepsilon,n},\phi_{\varepsilon,n})$)
such that
\begin{equation}
\Xi_n \in C^1([0,t_n]; V_n), \quad \quad \phi_n \in C^1([0,t_n]; V_n),
\end{equation}
and, for every $v \in V_n$ and for every $t \in [0,t_n]$, solving the approximating problem $(P_{\varepsilon,n})$:
\begin{equation} \nonumber 
(\partial_t[\Xi_n(t) + (\ell-\alpha) \phi_n(t)] -k \Delta \Xi_n(t) + k \alpha \Delta \phi_n(t) + A_{\varepsilon} \Xi_n(t), v )_H 
\end{equation}
\begin{equation} \label{1}
= (f_{\varepsilon}(t)-k \Delta \Xi^{*}, v)_H, 
\end{equation}
\begin{equation} \nonumber
(\partial_t \phi_n (t)- \nu \Delta \phi_n (t)+ \chi_{\varepsilon}(\phi_n(t)) + \pi(\phi_n(t)), v)_H 
\end{equation}
\begin{equation} \label{2}
= (\gamma [\Xi_n(t) - \alpha \phi_n(t) + \Xi^{*}], v)_H, 
\end{equation}
\begin{equation}
\partial_{N}\Xi_n = 0, \quad \quad \partial_{N}\phi_n = 0 \quad \textrm{on $\Sigma$},
\end{equation}
\begin{equation} \label{333}
\Xi_n(0) = \Xi_{0,n}, \quad \quad \phi_n(0) = \phi_{0,n} \quad \textrm{in $\Omega$}.
\end{equation}
This is a Cauchy problem for a system of nonlinear ordinary differential equations. In the next section we will show by a 
change of variable that this system admits a local solution $(\Xi_n,\phi_n)$, which is of the form 
\begin{equation} \label{liscio1}
\phi_n(t)= \sum^{n}_{i=1} a_{in}(t) v_i, 
\end{equation}
\begin{equation} \label{liscio2}
\Xi_n(t) = \sum^{n}_{i=1} b_{in}(t) v_i, 
\end{equation}
for some $a_{in} \in C^1([0,t_n])$ and $b_{in} \in C^1([0,t_n])$.
\paragraph{Remark 3.2}
We point out that
\begin{equation} \label{betacappuccioinl1}
\int_{\Omega}  \widehat{\chi}_{\varepsilon}(\phi_{0,n}) \leq C + \frac{1}{2\varepsilon} \| \phi_0- \phi_{0,n} \|_H ( \|\phi_0 \|_H + \|\phi_{0,n}\|_H),
\end{equation}
where
\begin{equation} \label{betacappuccioinl12}
C =  \| \widehat{\chi}(\phi_0) \|_{L^1(\Omega)}.
\end{equation}
Indeed, for every $\varepsilon \in (0,1]$, thanks to the property \eqref{questaservepoi} of $\widehat{\chi}_{\varepsilon}$,  we  have that
\begin{equation} 
0 \leq \widehat{\chi}_{\varepsilon}(\phi_0) \leq  \widehat{\chi}(\phi_0).
\end{equation}
Since $\widehat{\chi}(\phi_0) \in L^1(\Omega)$ (see \eqref{betal10}), we obtain that 
\begin{equation} 
\int_{\Omega} \widehat{\chi}_{\varepsilon}(\phi_0) \leq C,
\end{equation}
where $C = \| \widehat{\chi}(\phi_0) \|_{L^1(\Omega)}$.
From \eqref{questaservepoimolto}, using the Lipschitz continuity of $\chi_{\varepsilon}$, we have
\begin{eqnarray} \nonumber 
\widehat{\chi}_{\varepsilon}(\phi_{0,n}) 
&\leq& \widehat{\chi}_{\varepsilon}(\phi_0) + \Big| \int_{\phi_0}^{\phi_{0,n}} \chi_{\varepsilon}(s) \ ds \Big| \\ \nonumber
&\leq& \widehat{\chi}_{\varepsilon}(\phi_0) + \frac{1}{\varepsilon} \int_{\phi_0}^{\phi_{0,n}} |s| \ ds  \\ \label{disequazionebeta1epsilon}
&\leq&   \widehat{\chi}_{\varepsilon}(\phi_0) + \frac{1}{2\varepsilon} |\phi_0- \phi_{0,n}| (|\phi_0| + |\phi_{0,n}|).  
\end{eqnarray}
By integrating \eqref{disequazionebeta1epsilon} over $\Omega$, we obtain  that
\begin{equation} \label{Qe}
\int_{\Omega} \widehat{\chi}_{\varepsilon}(\phi_{0,n}) \leq Q_{\varepsilon}(n),
\end{equation}
where
\begin{equation} \nonumber
Q_{\varepsilon}(n) = C + \frac{1}{2\varepsilon} \| \phi_0- \phi_{0,n} \|_H ( \|\phi_0 \|_H + \|\phi_{0,n}\|_H).
\end{equation}

\begin{Rem4}
\emph{Thanks to \eqref{liscio2} and the Lipschitz continuity of $A_{\varepsilon}$, we obtain that}
\begin{equation} \label{Aepsilon}
A_{\varepsilon}(\Xi_n) \in C^{0,1}([0,t_n];H).
\end{equation}
\emph{Indeed,} $\| v_i \|_H \leq \| v_i \|_V  =1 $, 
\emph{for all} $i \in \mathbb{N}$. \emph{Then we choose} $t,t' \in [0, t_n]$ \emph{and we have the following inequality:}
\begin{equation} \nonumber
\begin{array}{ll}
\| A_{\varepsilon}(\Xi_n(t)) - A_{\varepsilon}(\Xi_n(t')) \|_H & = \| A_{\varepsilon}(\sum^{n}_{i=1} b_{in}(t) v_i) - A_{\varepsilon}(\sum^{n}_{i=1} b_{in}(t') v_i) \|_H \\
& \leq \frac{1}{\varepsilon}  \| \sum^{n}_{i=1} (b_{in}(t) -  b_{in}(t')) v_i \|_H \\
& \leq \frac{1}{\varepsilon}  \sum^{n}_{i=1} |b_{in}(t) -  b_{in}(t')| \, \| v_i \|_H \\
& = \frac{1}{\varepsilon}  \sum^{n}_{i=1} |b_{in}(t) -  b_{in}(t')|. 
\end{array}
\end{equation}
\emph{Since} $b_{in}$ \emph{are continuous, we obtain \eqref{Aepsilon}.}
\end{Rem4}

\paragraph{Existence of a local solution for $(P_{\varepsilon,n})$.}
In order to prove the existence of a local solution $(\Xi_n,\phi_n)$ for the approximating problem $(P_{\varepsilon,n})$, we make a change of variable. We set
\begin{equation}
\mu _n = \Xi_n + (\ell-\alpha) \phi_n, \quad \quad 
\mu _{0,n} = \Xi_{0,n} + (\ell-\alpha) \phi_{0,n},
\end{equation}
and we prove that there exists a local solution $(\mu _n,\phi_n)$ of the problem 
\begin{equation}
\begin{array}{ll}
(\partial_t \mu _n -k \Delta \mu _n +k\ell \Delta \phi_n
+ A_{\varepsilon} (\mu _n - (\ell-\alpha) \phi_n),v)_H = (f_{\varepsilon}-k \Delta \Xi^{*},v)_H, \\
(\partial_t \phi_n - \nu \Delta \phi_n + \chi_{\varepsilon}(\phi_n) + \pi(\phi_n),v)_H = (\gamma [\mu _n - \ell \phi_n + \Xi^{*}],v)_H, \\
\phi_n(0)= \phi_{0,n}, \quad \quad \mu _n(0)= \mu _{0,n},
\end{array}
\end{equation}
whenever $v \in V_n$.
Re-arranging the above system in explicit form, we have
\begin{equation} \label{esistenzalocale1111}
\begin{array}{ll}
(\partial_t \mu _n,v)_H = (k \Delta \mu _n - k\ell \Delta \phi_n
- A_{\varepsilon} (\mu _n - (\ell-\alpha) \phi_n) + f_{\varepsilon}-k \Delta \Xi^{*},v)_H, \\
(\partial_t \phi_n,v)_H = (\nu \Delta \phi_n - \chi_{\varepsilon}(\phi_n) - \pi(\phi_n) + \gamma [\mu _n - \ell \phi_n + \Xi^{*}],v)_H, \\
\phi_n(0)= \phi_{0,n}, \quad \quad \mu _n(0)= \mu _{0,n},
\end{array}
\end{equation}
whenever $v \in V_n$.
Thanks to the initial hypotheses \eqref{parametri}--\eqref{0}, \eqref{beta}--\eqref{pi}
and to the regularity of $A_{\varepsilon}$ shown in \eqref{Aepsilon}, the right-hand side of \eqref{esistenzalocale1111} is a Lipschitz continuous function from $[0,t_n]$ to $\mathbb{R}^n$.
Consequently, there exists a local solution for the approximating problem $(P_{\varepsilon,n})$.

\subsection{Global a priori estimates}
In this section we obtain four a priori estimates inferred from the main equations of the approximating problem $(P_{\varepsilon,n})$ (see \eqref{1}--\eqref{333}).

In the remainder of the paper we often use the H$\ddot{o}$lder inequality and to the elementary Young inequalities in performing our a priori estimates.
In particular, let us recall that, for every $a,b >0$, $\alpha \in (0,1)$ and $\delta > 0$, we have that
\begin{equation} \label{dis1}
ab \leq \alpha a^{\frac{1}{\alpha}} + (1 - \alpha)b^{\frac{1}{1- \alpha}},
\end{equation}
\begin{equation} \label{dis2}
ab \leq \delta a^2 + \frac{1}{4\delta}b^2.
\end{equation}
In the following the small-case symbol $c$ stands for different constants which depend only on $\Omega$, on the final time $T$,
on the shape of the nonlinearities and on the constants and the norms of the functions involved in the assumptions of our statements.

\paragraph{First a priori estimate.}
We add $\nu \phi_n $ to both sides of \eqref{2} and we test \eqref{1} by $\Xi_n$ and \eqref{2} by $\partial_t \phi_n$, respectively.
Then we sum up and integrate over $Q_t$, $t \in (0,T]$. We obtain that
\begin{equation} \label{laprimainassoluto}
\begin{array}{ll} 
\frac{1}{2} \int_{\Omega} |\Xi_n(t)|^2 + (\ell-\alpha) \int_{Q_t} \partial_t\phi_n\Xi_n + k \int_{Q_t} |\nabla \Xi_n|^2 
- k \alpha \int_{Q_t}\nabla \phi_n \cdot \nabla \Xi_n + \int_{Q_t} A_{\varepsilon} \Xi_n\Xi_n
\\
+ \int_{Q_t}|\partial_t \phi_n|^2 + \frac{\nu}{2} \int_{\Omega}  |\phi_n(t)|^2 
+ \frac{\nu}{2} \int_{\Omega}  |\nabla \phi_n(t)|^2 
+\int_{Q_t} \partial_t \widehat{\chi_{\varepsilon}}(\phi_n)
\\ 
=  \frac{1}{2} \int_{\Omega} |\Xi_{0,n}|^2 + \frac{\nu}{2} \int_{\Omega} |\phi_{0,n}|^2 + \frac{\nu}{2} \int_{\Omega}  |\nabla \phi_{0,n}|^2 +
\int_{Q_t} (f_{\varepsilon}-k \Delta \Xi^{*})\Xi_n 
\\
+ \int_{Q_t}[\gamma \Xi_n + (\nu-  \alpha\gamma) \phi_n + \gamma\Xi^{*}]\partial_t \phi_n
- \int_{Q_t}\pi(\phi_n)\partial_t \phi_n. 
\end{array}
\end{equation}	
To estimate the last integral on the right-hand side of \eqref{laprimainassoluto}, we observe that $\pi$ is a Lipschitz continuous function with Lipschitz constant $C_\pi$.
Consequently we have that
\begin{equation} \label{prima}
\begin{array}{ll}
|\pi(\phi_n)| & \leq  |\pi(\phi_n) - \pi(0)| + |\pi(0)| \\
& \leq C_{\pi}|\phi_n|+ |\pi(0)| \\
& \leq C_1(|\phi_n| +1),
\end{array}
\end{equation}
where $C_1 = \max{ \{ C_{\pi} ; |\pi(0)|  \} }$. Due to \eqref{dis2} and \eqref{prima}, we obtain that
\begin{equation} \label{pipi}
\begin{array}{ll}
- \int_{Q_t}\pi(\phi_n)\partial_t \phi_n &\leq \int_{Q_t} |\pi(\phi_n)\partial_t \phi_n| \\
& \leq \int_{Q_t} C_1(|\phi_n| +1) |\partial_t \phi_n|  \\
& \leq \frac{1}{8}\int_{Q_t} |\partial_t \phi_n  |^2 + 2C_1^2 \int_{Q_t} (|\phi_n| +1)^2 \\
& = \frac{1}{8}\int_{Q_t} |\partial_t \phi_n  |^2 + 4C_1^2 \int_{Q_t} |\phi_n|^2 + c.
\end{array}
\end{equation}
Now, we recall that $A_{\varepsilon}$ is a maximal monotone operator and $A_{\varepsilon}(0)=0$. Hence we have that
\begin{equation} \label{pezzoutile1}
\int_{Q_t} A_{\varepsilon} \Xi_n\Xi_n \geq 0. 
\end{equation}
Using \eqref{pipi}--\eqref{pezzoutile1}, from \eqref{laprimainassoluto} we obtain that
\begin{equation} \nonumber
\frac{1}{2} \int_{\Omega}  |\Xi_n (t)|^2 + k \int_{Q_t} | \nabla \Xi_n |^2 
+ \int_{Q_t}|\partial_t \phi_n|^2 + \frac{\nu}{2} \int_{\Omega}  |\phi_n (t)|^2   
+ \frac{\nu}{2} \int_{\Omega} |\nabla \phi_n(t)|^2 
+\int_{\Omega} \widehat{\chi_{\varepsilon}}(\phi_n(t)) 
\end{equation}
\begin{equation} \nonumber
\leq c + \frac{1}{2} \int_{\Omega}  |\Xi_{0,n}|^2 + \frac{\nu}{2} \int_{\Omega}  |\phi_{0,n}|^2 +  
\frac{\nu}{2} \int_{\Omega} |\nabla \phi_{0,n}|^2 + \int_{\Omega} \widehat{\chi_{\varepsilon}}(\phi_{0,n})
\end{equation}
\begin{equation} \nonumber
- (\ell-\alpha) \int_{Q_t} \partial_t\phi_n\Xi_n 
+ k \alpha \int_{Q_t} \nabla \phi_n \cdot \nabla \Xi_n 
+ \frac{1}{8}\int_{Q_t} |\partial_t \phi_n  |^2 + 4C_1^2 \int_{Q_t} |\phi_n|^2
\end{equation}
\begin{equation} \label{riferimento} 
+ \int_{Q_t} (f_{\varepsilon}- k \Delta \Xi^{*})\Xi_n 
+ \int_{Q_t}[\gamma \Xi_n + (\nu-  \alpha\gamma) \phi_n + \gamma\Xi^{*}]\partial_t \phi_n.
\end{equation}
We notice that the convergence provided by \eqref{datiinizialifaedo} assures that $\Xi_{0,n}$ and $\phi_{0,n}$ are bounded in $V$.
Consequently, thanks to \eqref{Qe}, the first four integrals on the right-hand side of \eqref{riferimento} are estimated as follows:
\begin{equation} \label{Questaqeullautile}
\frac{1}{2} \int_{\Omega}  |\Xi_{0,n}|^2 + \frac{\nu}{2} \int_{\Omega}  |\phi_{0,n}|^2 +  
\frac{\nu}{2} \int_{\Omega} |\nabla \phi_{0,n}|^2 + \int_{\Omega} \widehat{\chi_{\varepsilon}}(\phi_{0,n}) \leq c + Q_{\varepsilon}(n).
\end{equation}
We also notice that
\begin{eqnarray} \label{stimasuigradienti}
k \alpha \int_{Q_t} \nabla \phi_n \cdot \nabla \Xi_n \nonumber
&=& \frac{k}{2} \Big( 2 \alpha \int_{Q_t} \nabla \phi_n \cdot \nabla \Xi_n \Big) \\ \nonumber
&\leq & \frac{k}{2} \int_{Q_t} |\nabla \Xi_n|^2 + \frac{k\alpha^2}{2} \int_{Q_t} |\nabla \phi_n|^2   \\ 
&= & \frac{k}{2} \int_{Q_t} |\nabla \Xi_n|^2 + \frac{k\alpha^2}{\nu} \int_{Q_t} \frac{\nu}{2} |\nabla \phi_n|^2.   
\end{eqnarray} 
We re-arrange the right-hand side of \eqref{riferimento} using \eqref{dis2}, \eqref{Questaqeullautile}
and \eqref{stimasuigradienti}. Then we have that
\begin{equation} \nonumber
\frac{1}{2} \| \Xi_n (t) \|_H^2 + k \int_{Q_t} | \nabla \Xi_n |^2 
+ \int_{Q_t}|\partial_t \phi_n|^2 + \frac{\nu}{2} \|\phi_n (t)\|^2_V 
+\int_{\Omega} \widehat{\chi_{\varepsilon}}(\phi_n(t)) 
\end{equation}
\begin{equation} \nonumber
\leq c + Q_{\varepsilon}(n) + 2(\ell-\alpha)^2 \int_{Q_t} |\Xi_n|^2 + \frac{1}{8} \int_{Q_t} |\partial_t\phi_n|^2 
+ \frac{k}{2} \int_{Q_t} |\nabla \Xi_n|^2 + \frac{k\alpha^2}{\nu} \int_{Q_t} \frac{\nu}{2} |\nabla \phi_n|^2
\end{equation}
\begin{equation} \nonumber
+ \frac{1}{8}\int_{Q_t} |\partial_t \phi_n  |^2 + 4C_1^2 \int_{Q_t} |\phi_n|^2 
+ 2 \int_{Q_t} |f_{\varepsilon}- k \Delta \Xi^{*}|^2 
+ \frac{1}{8} \int_{Q_t} |\Xi_n|^2
\end{equation}
\begin{equation} \label{eccola111111}
+ 2 \int_{Q_t} |\gamma \Xi_n + (\nu-  \alpha\gamma) \phi_n + \gamma\Xi^{*}|^2  + \frac{1}{8} \int_{Q_t} |\partial_t \phi_n|^2.  
\end{equation}
According to \eqref{convergenzadif}, $f_{\varepsilon}$ is bounded in $L^2(0,T; H)$ uniformly with respect to $\varepsilon$.
Consequently, due to \eqref{f}--\eqref{eta}, the seventh integral on the right-hand side of \eqref{eccola111111} is under control 
and similarly the third addendum in the ninth integral on the right-hand side. Then we infer that
\begin{equation} \nonumber
\frac{1}{2} \| \Xi_n (t) \|_H^2 + k \int_{Q_t} | \nabla \Xi_n |^2 
+ \int_{Q_t}|\partial_t \phi_n|^2 + \frac{\nu}{2} \|\phi_n (t)\|^2_V 
+\int_{\Omega} \widehat{\chi_{\varepsilon}}(\phi_n(t)) 
\end{equation}
\begin{equation} \nonumber
\leq c + Q_{\varepsilon}(n) + \Big[ 2(\ell-\alpha)^2 + \frac{1}{8} \Big] \int_{Q_t} |\Xi_n|^2 
+ \frac{1}{2} \int_{Q_t} |\partial_t\phi_n|^2 
+ \frac{k\alpha^2}{\nu} \int_{Q_t} \frac{\nu}{2} |\nabla \phi_n|^2  
\end{equation}
\begin{equation} \label{eccola1111112}
+ 4C_1^2 \int_{Q_t} |\phi_n|^2 
+ 8\gamma^2 \int_{Q_t} |\Xi_n|^2 + 8(\nu-  \alpha\gamma)^2 \int_{Q_t} |\phi_n|^2 .
\end{equation}
Now, we recollect the constants in \eqref{eccola1111112} and obtain that 
\begin{equation} \nonumber
\frac{1}{2} \| \Xi_n (t) \|_H^2 + k \int_{Q_t} | \nabla \Xi_n |^2 
+ \int_{Q_t}|\partial_t \phi_n|^2 + \frac{\nu}{2} \|\phi_n (t)\|^2_V 
+\int_{\Omega} \widehat{\chi_{\varepsilon}}(\phi_n(t)) 
\end{equation}
\begin{equation} \label{eccola1111113}
\leq c + Q_{\varepsilon}(n) + C_2 \frac{1}{2} \int_0^t \| \Xi_n(s) \|_H^2 \ ds
+ C_3 \frac{\nu}{2} \int_0^t \| \nabla \phi_n(s) \|_H^2   \ ds
+ C_4 \frac{\nu}{2} \int_0^t \|\phi_n(s) \|_H^2 \ ds,
\end{equation}
where
\begin{equation} \nonumber
C_2= 2 [2(\ell-\alpha)^2 + \frac{1}{8} + 8\gamma^2], \quad C_3= \frac{k\alpha^2}{\nu}, \quad C_4= \frac{2[4C_1^2 + 8(\nu-  \alpha\gamma)^2  ]}{\nu}.
\end{equation}
Consequently, from \eqref{eccola1111113} we have that
\begin{equation} \nonumber
\frac{1}{2} \| \Xi_n (t) \|_H^2 + k \int_{Q_t} | \nabla \Xi_n |^2 
+ \int_{Q_t}|\partial_t \phi_n|^2 + \frac{\nu}{2} \|\phi_n (t)\|^2_V 
+\int_{\Omega} \widehat{\chi_{\varepsilon}}(\phi_n(t)) 
\end{equation}
\begin{equation} \label{cara2345678900}
\leq c + Q_{\varepsilon}(n) + C_5 \Bigg( \frac{1}{2} \int^t_0 \|\Xi_n (s) \|_H^2 \ ds + \frac{\nu}{2} \int^t_0 \|\phi_n(s)\|_V^2 \ ds \Bigg),
\end{equation}
where
\begin{equation} \nonumber
C_5 = \max{(C_2, \ C_3, \ C_4)}.
\end{equation} 
Then, from \eqref{cara2345678900} we conclude that
\begin{equation} \nonumber
\frac{1}{2} \| \Xi_n (t) \|_H^2 + k \int_{Q_t} | \nabla \Xi_n |^2 
+ \int_{Q_t}|\partial_t \phi_n|^2 + \frac{\nu}{2} \|\phi_n (t)\|^2_V 
+\int_{\Omega} \widehat{\chi_{\varepsilon}}(\phi_n(t)) 
\end{equation}
\begin{equation} \label{cara23456789}
\leq c_{\varepsilon} \Bigg( 1 + \frac{1}{2} \int^t_0 \|\Xi_n (s) \|_H^2 \ ds + \frac{\nu}{2} \int^t_0 \|\phi_n(s)\|_V^2 \ ds \Bigg).
\end{equation}
Now, we apply the Gronwall lemma to \eqref{cara23456789} and infer that
\begin{equation} \label{eccola1111114}
\frac{1}{2} \| \Xi_n (t) \|_H^2 + k \int_{Q_t} | \nabla \Xi_n |^2 
+ \int_{Q_t}|\partial_t \phi_n|^2 + \frac{\nu}{2} \|\phi_n (t)\|^2_V 
+\int_{\Omega} \widehat{\chi_{\varepsilon}}(\phi_n(t)) \leq c_{\varepsilon}.
\end{equation}
As \eqref{eccola1111114} holds true for any $t \in [0,t_n)$, we conclude that
\begin{equation} \label{1a}
\| \phi_n \|_{H^{1}(0,t_n; H) \cap L^{\infty}(0,t_n; V)} \leq c_{\varepsilon}, 
\end{equation}
\begin{equation} \label{1b}
\| \Xi_n \|_{L^{\infty}(0,t_n; H) \cap L^{2}(0,t_n; V)} \leq c_{\varepsilon}, 
\end{equation}
\begin{equation} \label{1c}
\| \widehat{\chi_{\varepsilon}}(\phi_n) \|_{L^{\infty}(0,t_n; L^1(\Omega))} \leq c_{\varepsilon}. 
\end{equation}

\paragraph{Second a priori estimate.}
First of all, we notice that $\pi(\phi_n)$ is bounded in $L^{2}(0,t_n; H)$ owing to \eqref{pi} and \eqref{1a}.
Thanks to \eqref{1a}--\eqref{1c}, we can rewrite \eqref{2} as
\begin{equation} \label{seconda2}
(- \nu \Delta \phi_n + \chi_{\varepsilon}(\phi_n),v)_H  = (g_1,v)_H, \quad \quad \textrm{for all $v \in V_n$,}
\end{equation}
with $\| g_1 \|_{L^{2}(0,t_n; H)} \leq c_{\varepsilon} $. 
The choice of the basis $v_i$ as in \eqref{baseortonormale} allows us to
test \eqref{seconda2} by $- \Delta \phi_n$. Integrating over $(0,t)$, we obtain that
\begin{equation} \label{siiiii1}
\nu \int_{Q_t}   | \Delta \phi_n |^2 + \int_{Q_t} \nabla \phi_n \cdot \nabla \chi_{\varepsilon}(\phi_n) = - \int_{Q_t} g_1 \Delta \phi_n.
\end{equation}
Using inequalities \eqref{dis1}--\eqref{dis2}, from \eqref{siiiii1} we have that
\begin{equation} \label{siiiii13456789}
\frac{\nu}{2} \int_{Q_t}   | \Delta \phi_n |^2 + \int_{Q_t} \chi'_{\varepsilon}(\phi_n)|\nabla \phi_n|^2  
\leq \frac{1}{2 \nu } \int_{Q_t} |g_1|^2. 
\end{equation}
Due to \eqref{1a} and the monotonicity of $\chi_{\varepsilon}$, from \eqref{siiiii13456789} we obtain that
\begin{equation} \label{2aaaaa}
\| \Delta \phi_n \|_{L^{2}(0,t; H)} \leq c_{\varepsilon}.
\end{equation}
We observe that \eqref{2aaaaa} holds true for any $t \in [0,t_n)$.
Then, using elliptic regularity, from \eqref{1a} and \eqref{2aaaaa} we infer that
\begin{equation} \label{2a}
\| \phi_n \|_{L^{2}(0,t_n; W)} \leq c_{\varepsilon}.
\end{equation}

\paragraph{Third a priori estimate.}
Thanks to the previous a priori estimates, from \eqref{1} it follows that
\begin{equation} \label{nuovarelazione}
(\partial_t \Xi_n  -k \Delta \Xi_n + A_{\varepsilon} \Xi_n,v)_H =  (g_2,v)_H \quad \quad \textrm{for all $v \in V_n$,}
\end{equation}
with $\| g_2 \|_{L^{2}(0,t_n; H)} \leq c_{\varepsilon} $. We test \eqref{nuovarelazione} by $\partial_t \Xi_n$ and integrate over $(0,t)$; we obtain that
\begin{equation} \label{aggiunta11111}
\int_{Q_t} |\partial_t \Xi_n|^2  + \frac{k}{2} \int_{\Omega} | \nabla \Xi_n(t)|^2 + \int_{Q_t} A_{\varepsilon} \Xi_n  \partial_t \Xi_n = 
\frac{k}{2} \int_{\Omega} | \nabla \Xi_{0,n}|^2 + \int_{Q_t} g_2  \partial_t \Xi_n.
\end{equation}
Then, using the property \eqref{Ae} of $A_{\varepsilon}$
and inequalities \eqref{dis1}--\eqref{dis2}, from \eqref{aggiunta11111} we infer that
\begin{equation} \nonumber
\int_{Q_t} |\partial_t \Xi_n|^2  + \frac{k}{2} \int_{\Omega} | \nabla \Xi_n (t)| 
\end{equation}
\begin{eqnarray} \nonumber
&\leq& \frac{k}{2} \int_{\Omega} | \nabla \Xi_{0,n}|^2  +  \int_{Q_t} |A_{\varepsilon} \Xi_n  \partial_t \Xi_n| + 2 \int_{Q_t} |g_2|^2 + \frac{1}{8} \int_{Q_t}  |\partial_t \Xi_n|^2 \\ \nonumber
&\leq & \frac{k}{2} \int_{\Omega} | \nabla \Xi_{0,n}|^2 +  2 \int_{Q_t} | A_{\varepsilon} \Xi_n |^2 + \frac{1}{8} \int_{Q_t} |\partial_t \Xi_n |^2 
+ 2 \int_{Q_t} |g_2|^2 + \frac{1}{8} \int_{Q_t}  |\partial_t \Xi_n|^2 \\ \nonumber
& = & \frac{k}{2} \int_{\Omega} | \nabla \Xi_{0,n}|^2 +  2 \int^t_0  \| A_{\varepsilon} \Xi_n(s) \|^2_H \ ds+ \frac{1}{4} \int_{Q_t} |\partial_t \Xi_n |^2 + 2 \int_{Q_t} |g_2|^2 \\ \nonumber
& \leq & \frac{k}{2} \int_{\Omega} | \nabla \Xi_{0,n}|^2 +  2 \int^t_0 [C (\| \Xi_n (s)\|_H + 1)]^2 \ ds+ \frac{1}{4} \int_{Q_t} |\partial_t \Xi_n |^2 + 2 \int_{Q_t} |g_2|^2  \\ \label{aggiunta111113}
& \leq & c + \frac{k}{2} \int_{\Omega} | \nabla \Xi_{0,n}|^2 + 4 C^2 \int^t_0 \| \Xi_n(s) \|^2_H \ ds + \frac{1}{2} \int_{Q_t} |\partial_t \Xi_n |^2 + 2 \int_{Q_t} |g_2|^2. 
\end{eqnarray} 
Due to \eqref{0}, the first integral on the right-hand side of \eqref{aggiunta111113} is under control. Then,
from \eqref{aggiunta111113} we infer that
\begin{equation} \label{carissima111111}
\frac{1}{2} \int_{Q_t} |\partial_t \Xi_n|^2  + \frac{k}{2} \int_{\Omega} | \nabla \Xi_n (t)| 
\leq c +  4C^2 \int^t_0 \| \Xi_n(s) \|^2_H \ ds + 2 \int_{Q_t} |g_2|^2.
\end{equation}
We observe that \eqref{carissima111111} holds true for any $t \in [0,t_n)$,
Then, due to the previous estimates \eqref{1a}--\eqref{1b}, we conclude that
\begin{equation} \label{3a}
\| \Xi_n \|_{H^{1}(0,t_n; H) \cap L^{\infty}(0,t_n; V)} \leq c_{\varepsilon}.
\end{equation}

\paragraph{Fourth a priori estimate.}
Due to the previous estimates \eqref{1a}--\eqref{1c}, \eqref{2a} and \eqref{3a}, by comparison in \eqref{nuovarelazione}, we infer that 
\begin{equation}
\| \Delta \Xi_n \|_{L^{2}(0,t_n; H)} \leq c_{\varepsilon}.
\end{equation}
Consequently,  we conclude that
\begin{equation} \label{4a}
\| \Xi_n \|_{L^{2}(0,t_n; W)} \leq c_{\varepsilon}.
\end{equation}

\paragraph{Summary of the a priori estimates.}
Since the constants appearing in the a priori estimates are all independent of $t_n$,
the local solution can be extended to a solution defined on the whole interval $[0,T]$, i.e.,
we can assume $t_n = T$ for any $n$. Hence, due to \eqref{1a}--\eqref{1c}, \eqref{2a}, \eqref{3a} and \eqref{4a}, we conclude that  
\begin{eqnarray} 
\label{1N}
\| \phi_n \|_{H^{1}(0,T; H) \cap L^{\infty}(0,T; V) \cap L^{2}(0,T; W)} \leq c_{\varepsilon}, \\
\label{2N}
\| \Xi_n \|_{H^{1}(0,T; H) \cap L^{\infty}(0,T; V) \cap L^{2}(0,T; W)} \leq c_{\varepsilon}.
\end{eqnarray}

\subsection{Passage to the limit as $n \rightarrow + \infty$}
Now, we let $n \rightarrow + \infty$ and show that 
the limit of some subsequences of solutions for $(P_{\varepsilon, n})$ (see \eqref{1}--\eqref{333}) yields a solution 
of $(P_{\varepsilon})$ (see \eqref{equazionedelsecondolimite1}--\eqref{equazionedelsecondolimite2}).
Estimates \eqref{1N}--\eqref{2N} for $\phi_n$ and $\Xi_n$ and the well-known weak or weak*
compactness results ensure the existence of a pair $(\phi_{\varepsilon},\Xi_{\varepsilon})$ such that,
at least for a subsequence, 
\begin{eqnarray} \label{necessaria1}
\phi_n \rightharpoonup \phi_{\varepsilon} & \textrm{in}  & H^1(0,T;H) \cap L^2(0,T;W), \\
\phi_n \rightharpoonup^{*} \phi_{\varepsilon} & \textrm{in}  & L^{\infty}(0,T;V),\\ 
\Xi_n \rightharpoonup \Xi_{\varepsilon} & \textrm{in}  & H^1(0,T;H) \cap L^2(0,T;W), \\ \label{necessaria2}
\Xi_n \rightharpoonup^{*} \Xi_{\varepsilon} & \textrm{in}  & L^{\infty}(0,T;V),
\end{eqnarray}
as $n \rightarrow + \infty$.
We notice that $W$, $V$, $H$ are Banach spaces and
$W \subset V \subset H$ with dense and compact embeddings.
Then, we are under the assumptions of \cite[Prop.~4, Sec.~8]{Simon} and this fact implies the following strong convergences:
\begin{eqnarray} \label{conv1}
\phi_n \rightarrow \phi_{\varepsilon} & \textrm{in}  & C^0([0,T];H) \cap L^2(0,T;V), \\ 
\Xi_n \rightarrow \Xi_{\varepsilon} & \textrm{in}  & C^0([0,T];H) \cap L^2(0,T;V), \label{conv2} 
\end{eqnarray}
as $n \rightarrow + \infty$. Since $\pi$, $A_{\varepsilon}$ and $\chi_{\varepsilon}$  are Lipschitz continuous, we infer that
\begin{equation} \label{correzione22222222}
\begin{array}{ll}
|\pi(\phi_n) - \pi(\phi_{\varepsilon}) | \leq C_{\pi} |\phi_n - \phi_{\varepsilon} |  & \textrm{a.e. in $Q$,} \\
\| A_{\varepsilon}\Xi_n - A_{\varepsilon}\Xi_{\varepsilon}  \|_H \leq \frac{1}{\varepsilon} \| \Xi_n - \Xi_{\varepsilon} \|_H & \textrm{a.e. in $[0,T]$,} \\
|\chi_{\varepsilon}(\phi_n) - \chi_{\varepsilon}(\phi_{\varepsilon}) | \leq \frac{1}{\varepsilon} |\phi_n - \phi_{\varepsilon} |  & \textrm{a.e. in $Q$.} \\
\end{array}
\end{equation}
Due to \eqref{correzione22222222}, we conclude that 
\begin{eqnarray}
\pi(\phi_n) \rightarrow \pi(\phi_{\varepsilon}) & \textrm{in}  & C^0([0,T];H), \\
A_{\varepsilon}\Xi_n \rightarrow A_{\varepsilon}\Xi_{\varepsilon} & \textrm{in}  & C^0([0,T];H), \\
\chi_{\varepsilon}(\phi_n) \rightarrow \chi_{\varepsilon}(\phi_{\varepsilon}) & \textrm{in}  & C^0([0,T];H), \label{necessaria3}
\end{eqnarray}
as $n \rightarrow + \infty$.
Now, we fix $k \leq n$ and we observe that, for every $v \in V_k$ and for every $t \in [0,T]$, 
the solution $(\Xi_n,\phi_n)$ of problem $(P_{\varepsilon,n})$ satisfies
\begin{equation} \nonumber
(\partial_t[\Xi_n(t) + (\ell-\alpha) \phi_n(t)] -k \Delta \Xi_n(t) + k \alpha \Delta \phi_n(t) + A_{\varepsilon} \Xi_n(t), v )_H 
\end{equation}
\begin{equation} \label{amico1}
= (f_{\varepsilon}(t)-k \Delta \Xi^{*}, v)_H, 
\end{equation}
\begin{equation} \nonumber
(\partial_t \phi_n (t)- \nu \Delta \phi_n (t)+ \chi_{\varepsilon}(\phi_n(t)) + \pi(\phi_n(t)), v)_H 
\end{equation}
\begin{equation} \label{amico2}
= (\gamma [\Xi_n(t) - \alpha \phi_n(t) + \Xi^{*}], v)_H.
\end{equation}
If $k$ is fixed and $n \rightarrow + \infty$, we have the convergence of every term of \eqref{amico1}--\eqref{amico2}
to the corresponding one with $\Xi_{\varepsilon}$, $\phi_{\varepsilon}$ whenever $v \in V_k$, i.e.,
\begin{equation}\nonumber
(\partial_t[\Xi_{\varepsilon}(t) + (\ell-\alpha) \phi_{\varepsilon}(t)] -k \Delta \Xi_{\varepsilon}(t) + k \alpha \Delta \phi_{\varepsilon}(t) + A_{\varepsilon} \Xi_{\varepsilon}(t), v )_H 
\end{equation}
\begin{equation} 
= (f_{\varepsilon}(t)-k \Delta \Xi^{*}, v)_H, 
\end{equation}
\begin{equation} \nonumber
(\partial_t \phi_{\varepsilon} (t)- \nu \Delta \phi_{\varepsilon} (t)+ \chi_{\varepsilon}(\phi_{\varepsilon}(t)) + \pi(\phi_{\varepsilon}(t)), v)_H 
\end{equation}
\begin{equation} 
= (\gamma [\Xi_{\varepsilon}(t) - \alpha \phi_{\varepsilon}(t) + \Xi^{*}], v)_H.
\end{equation}
As $k$ is arbitrary, the limit equalities hold true
for every $v \in \bigcup^{\infty}_{k=1} V_k$, which is dense in $V$. Then the limit equalities actually hold for every $v \in V$, i.e.,
\begin{equation}
\partial_t(\Xi_{\varepsilon} + (\ell-\alpha) \phi_{\varepsilon}) -k \Delta \Xi_{\varepsilon} + k \alpha \Delta \phi_{\varepsilon} + A_{\varepsilon} \Xi_{\varepsilon} = f_{\varepsilon}-k \Delta \Xi^{*} \ \ \textrm{a.e. in $Q$,}
\end{equation}
\begin{equation} \label{utilenelsecondolimite}
\partial_t \phi_{\varepsilon} - \nu \Delta \phi_{\varepsilon} + \chi_{\varepsilon}(\phi_{\varepsilon}) + \pi(\phi_{\varepsilon}) = \gamma (\Xi_{\varepsilon} - \alpha \phi_{\varepsilon} + \Xi^{*})  \ \ \textrm{a.e. in $Q$.}
\end{equation}
Now, we prove the convergence of the initial data. We recall that
\begin{equation}
\Xi_{0,n} = P_{V_n} \Xi_0, 
\quad \quad \quad \quad
\phi_{0,n} = P_{V_n} \phi_0. 
\end{equation}
If $\varepsilon$ is fixed, then
\begin{equation} \label{questopountodelicato1}
\lim_{n \rightarrow + \infty} \Xi_{0,n} =  \Xi_0 \quad \textrm{in $V$}, 
\end{equation} 
\begin{equation} \label{questopountodelicato2}
\lim_{n \rightarrow + \infty} \phi_{0,n} = \phi_0 \quad \textrm{in $V$},
\end{equation}
and then also in $H$.
These observations and \eqref{conv1}--\eqref{conv2}
show that the weak limit of some subsequences of solutions for $(P_{\varepsilon, n})$ (see \eqref{1}--\eqref{333})
yields a solution for $(P_{\varepsilon})$ (see \eqref{equazionedelsecondolimite1}--\eqref{equazionedelsecondolimite2}).
We also notice that taking the limit as $n \rightarrow + \infty$ in \eqref{Qe} entails that $ Q_{\varepsilon}(n) \rightarrow C$, with
\begin{equation} 	\label{amico33333}
\int_{\Omega} \widehat{\chi}_{\varepsilon}(\phi_0) \leq C.
\end{equation}
Then, after the first passage to the limit, we conclude that 
estimates \eqref{1N}--\eqref{2N} still hold for the limiting functions with constants independent of $\varepsilon$, i.e.,
\begin{eqnarray} 
\label{1Nuovissima}
\| \phi_{\varepsilon} \|_{H^{1}(0,T; H) \cap L^{\infty}(0,T; V) \cap L^{2}(0,T; W)} \leq c, \\
\label{2Nuovissima}
\| \Xi_{\varepsilon} \|_{H^{1}(0,T; H) \cap L^{\infty}(0,T; V) \cap L^{2}(0,T; W)} \leq c.
\end{eqnarray}


\subsection{Passage to the limit as $\varepsilon \searrow 0$}
Now, we let $\varepsilon  \searrow 0$ and show that the limit of some subsequences of solutions for $(P_{\varepsilon})$
(see \eqref{equazionedelsecondolimite1}--\eqref{equazionedelsecondolimite2}) tends to a solution 
of the initial problem $(P)$ (see \eqref{iniziale1}--\eqref{iniziale5}). First of all, due to \eqref{necessaria1}--\eqref{conv2}, \eqref{necessaria3} 
and \eqref{amico33333}, we have that the constants in \eqref{1Nuovissima}--\eqref{2Nuovissima}
do not depend on $\varepsilon$. Moreover, thanks to \eqref{1Nuovissima}--\eqref{2Nuovissima}, by comparison in \eqref{utilenelsecondolimite}, we infer that
\begin{equation} \label{3Ne}
\| \chi_{\varepsilon}(\phi_{\varepsilon}) \|_{L^2(Q)} \leq c. 
\end{equation} 
The well-known weak or wea$\emph{k}^{*}$ compactness results and 
the useful theorem \cite[Prop.~4, Sec.~8]{Simon} ensure the existence of a pair $(\phi,\Xi)$ such that,
at least for a subsequence, 
\begin{eqnarray} 
\phi_{\varepsilon} \rightharpoonup^{*} \phi & \textrm{in  $ H^1(0,T;H) \cap L^{\infty}(0,T; V) \cap L^2(0,T;W)$}, &\\
\Xi_{\varepsilon} \rightharpoonup^{*} \Xi       & \textrm{in  $ H^1(0,T;H) \cap L^{\infty}(0,T; V) \cap L^2(0,T;W)$},  &\\ \label{convutilissimaadesso}
\phi_{\varepsilon} \rightarrow \phi         & \textrm{ in $ C^0([0,T];H) \cap L^2(0,T;V)$}, & \\ \label{convutilissimaadesso1111}
\Xi_{\varepsilon}    \rightarrow \Xi            & \textrm{  in $C^0([0,T];H) \cap L^2(0,T;V)$},&
\end{eqnarray} 
as $\varepsilon \searrow 0$.
Now, we observe that \eqref{convutilissimaadesso} implies that
\begin{equation} \label{amico44}
\phi_{\varepsilon} \rightarrow \phi  \quad \quad \textrm{in $L^2(0,T;H) \equiv L^2(Q)$}  
\end{equation}
as $\varepsilon \searrow 0$.
We set $\xi_{\varepsilon} = \chi_{\varepsilon}(\phi_{\varepsilon})$ and remark that
\begin{equation} 
\| \xi_{\varepsilon} \|_{L^{2}(Q)} = \| \chi_{\varepsilon}(\phi_{\varepsilon}) \|_{L^{2}(Q)} \leq c.
\end{equation}
Thus, we may suppose that, as $\varepsilon \searrow 0$,
at least for a subsequence,
\begin{equation} \label{deboleepsilon} 
\xi_{\varepsilon} \rightharpoonup \xi \quad \quad \textrm{in} \ L^{2}(Q),
\end{equation}
for some $\xi \in L^{2}(Q)$. 
Now, we introduce the operator $\mathcal{B}_{\varepsilon}$ induced by $\chi_{\varepsilon}$ on $L^2(Q)$ in the following way:
\begin{equation} 
\mathcal{B}_{\varepsilon}: L^2(Q) \longrightarrow L^2(Q)
\end{equation}
\begin{equation} 
\xi_{\varepsilon}  \in \mathcal{B}_{\varepsilon}(\phi_{\varepsilon}) \Longleftrightarrow \xi_{\varepsilon}(x,t) \in \chi_{\varepsilon}(\phi_{\varepsilon}(x,t)) \quad \textrm{for a.e. $(x,t) \in Q$.}
\end{equation}
Due to  \eqref{amico44} and \eqref{deboleepsilon}, we have that
\begin{equation} \label{fine11}
\left\{ \begin{array}{ll}
\mathcal{B}_{\varepsilon}(\phi_{\varepsilon}) \rightharpoonup \xi & \textrm{in $L^2(Q)$,} \\
\phi_{\varepsilon} \rightarrow \phi & \textrm{in $L^2(Q)$,}
\end{array}
\right.
\end{equation}
\begin{equation} \label{fine12}
\limsup_{\varepsilon \searrow 0} \int_Q \xi_{\varepsilon} \phi_{\varepsilon} = \int_Q \xi \phi.
\end{equation}
Thanks to \eqref{fine11}--\eqref{fine12} and to the useful results proved in \cite[Prop.~2.2, p.~38]{Barbu}, we conclude that
\begin{equation} \label{proprioquestaqui}
\xi  \in \mathcal{B}(\phi) \quad \textrm{in $L^2(Q)$},
\end{equation}
where $\mathcal{B}$ is defined by \eqref{betagrande1}--\eqref{betagrande2}.
This is equivalent to say that
\begin{equation} 
\xi  \in \chi(\phi) \quad \textrm{a.e. in $Q$}. 
\end{equation}
Moreover, we pass to the limit in $A_{\varepsilon}$ by repeating
the previous arguments and conclude that
\begin{equation} 
\zeta \in \mathcal{A}(\Xi) \quad \textrm{in $L^2(0,T;H)$}, 
\end{equation}
with obvious definition for $\mathcal{A}$ (see \eqref{Agrande1}--\eqref{Agrande2}), and this is equivalent to say that
\begin{equation} 
\zeta \in A(\Xi) \quad \textrm{a.e. in $[0,T]$.} 
\end{equation}

\paragraph{Conclusion of the proof.}
Thanks to the previous steps, 
we conclude that, as $\varepsilon \searrow 0$, the limit of some subsequences of solutions $(\Xi_{\varepsilon},\phi_{\varepsilon})$ to $(P_{\varepsilon})$
(see \eqref{equazionedelsecondolimite1}--\eqref{equazionedelsecondolimite2})
yields a solution $(\Xi,\phi)$ of the initial boundary value problem $(P)$, i.e.,
\begin{equation} 
\partial_t(\Xi + (\ell-\alpha) \phi) -k \Delta \Xi + k \alpha \Delta \phi + \zeta = f-k \Delta \Xi^{*} \ \ \textrm{a.e. in $Q$,} 
\end{equation}
\begin{equation} 
\partial_t \phi - \nu \Delta \phi + \xi + \pi(\phi) = \gamma (\Xi - \alpha \phi + \Xi^{*})  \ \ \textrm{a.e. in $Q$,}
\end{equation}
\begin{equation} 
\zeta(t) \in A(\Xi(t)) \ \ \textrm{for a.e. $ t \in (0,T)$,}  
\end{equation}
\begin{equation} 
\xi \in \chi(\phi) \ \ \textrm{a.e. in $Q$,}
\end{equation}
\begin{equation}
\partial_{N}\Xi = 0, \quad \quad \partial_{N}\phi = 0 \quad \textrm{on $\Sigma$},
\end{equation}
\begin{equation}
\Xi(0) = \Xi_0, \quad \quad \phi(0) = \phi_0 \quad \textrm{in $\Omega$.}
\end{equation}
We notice that the homogeneous Neumann boundary conditions for both $\Xi$ and $\phi$ 
follow from \eqref{regloarita}, due to the definition of $W$ (see \eqref{W}).

\setcounter{equation}{0}

\section{Proof of the continuous dependence theorem}
This section is devoted to the proof of Theorem \ref{Teo3}.

Assume $\alpha = \ell$. If $f_i$, $\Xi^{*}_i$, $\Xi_{0_i}$, $\phi_{0_i}$, $i=1,2$, 
are given as in \eqref{f}--\eqref{0} and $(\phi_i, \Xi_i)$, $i=1,2$,  are the corresponding solutions, 
we can write problem \eqref{iniziale1}--\eqref{iniziale5} for both $(\phi_i, \Xi_i)$, $i=1,2$, obtaining
\begin{equation} \label{uni1}
\partial_t\Xi_i -k \Delta \Xi_i + k \ell \Delta \phi_i + \zeta_i = f_i-k \Delta \Xi^{*}_i \ \ \textrm{a.e. in $Q$,}
\end{equation}
\begin{equation} \label{uni2}
\partial_t \phi_i - \nu \Delta \phi_i + \xi_i + \pi(\phi_i) = \gamma (\Xi_i - \ell \phi_i + \Xi^{*}_i) \ \ \textrm{a.e. in $Q$,}
\end{equation}
\begin{equation} 
\zeta_i(t) \in A(\Xi_i(t)) \ \ \textrm{for a.e. $t \in (0,T)$,}  
\end{equation}
\begin{equation} 
\xi_i \in \chi(\phi_i) \ \ \textrm{a.e. in $Q$,}
\end{equation}
\begin{equation}
\partial_{N}\Xi_i = 0, \quad \quad \partial_{N}\phi_i = 0 \quad \textrm{on $\Sigma$},
\end{equation}
\begin{equation} \label{ultimafinale}
\Xi_i(0) = \Xi_{0_i}, \quad \quad \phi_i(0) = \phi_{0_i}.
\end{equation}
First of all, we set
\begin{eqnarray}
\phi=\phi_1 - \phi_2, & \Xi= \Xi_1 - \Xi_2, \\
f= f_1 - f_2,  & \Xi^* = \Xi^*_1 - \Xi^*_2, \\ 
\phi_0= \phi_{0_1} - \phi_{0_2}, & \Xi_0= \Xi_{0_1} - \Xi_{0_2}.
\end{eqnarray}
We write \eqref{uni1} for both $(\phi_1, \Xi_1)$ and $(\phi_2, \Xi_2)$ and we take the difference.
We obtain that
\begin{equation} \label{prec1}
\partial_t \Xi -k \Delta \Xi + k \ell \Delta \phi + \zeta_1 - \zeta_2 = f-k \Delta \Xi^{*}. 
\end{equation} 
We write \eqref{uni2} for both $(\phi_1, \Xi_1)$ and $(\phi_2, \Xi_2)$ and we take the difference.
We obtain that
\begin{equation} \label{prec2}
\partial_t \phi - \nu \Delta \phi + \xi_1 - \xi_2 + \pi(\phi_1) - \pi(\phi_2) = \gamma (\Xi - \ell \phi + \Xi^{*}). 
\end{equation}
We multiply \eqref{prec1} by $\Xi$ and \eqref{prec2} by $\frac{kl^2}{\nu}\phi$.
Then we  sum up and integrate over $Q_t$, $t \in (0,T]$. We have that
\begin{equation} \nonumber
\frac{1}{2} \int_{\Omega} |\Xi(t)|^2 + \frac{k\ell^2}{2\nu} \int_{\Omega} |\phi(t)|^2 + k \int_{Q_t}  (|\nabla \Xi|^2 -\ell \nabla\phi\nabla\Xi + \ell^2 |\nabla\phi|^2)
\end{equation}
\begin{equation} \nonumber
+ \int_{Q_t} (\zeta_1 - \zeta_2) (\Xi_1 - \Xi_2) 
+ \frac{k\ell^2}{\nu} \int_{Q_t} (\xi_1 - \xi_2) (\phi_1 - \phi_2) 
\end{equation}
\begin{equation} \nonumber 
= \frac{1}{2} \| \Xi_0 \|^2_H + \frac{k\ell^2}{2\nu} \| \phi_0 \|^2_H 
- \frac{k\ell^2}{\nu} \int_{Q_t} [\pi(\phi_1) - \pi(\phi_2)](\phi_1 - \phi_2) 
\end{equation}
\begin{equation} \label{tutta1}
+ \int_{Q_t} (f-k \Delta \Xi^{*})\Xi 
+ \frac{\gamma k\ell^2}{\nu} \int_{Q_t} \Xi\phi-\frac{\gamma k\ell^3}{\nu} \int_{Q_t} |\phi|^2 + \frac{\gamma k\ell^2}{\nu} \int_{Q_t}\Xi^{*}\phi. 
\end{equation}
Since $A$ and $\chi$ are  maximal monotone, we have that
\begin{equation} \label{stimissimaamica1}
\int_{Q_t} (\zeta_1 - \zeta_2) (\Xi_1 - \Xi_2) \geq 0, 
\end{equation} 
\begin{equation} \label{stimissimaamica2}
\int_{Q_t} (\xi_1 - \xi_2) (\phi_1 - \phi_2) \geq 0.
\end{equation} 
Moreover, thanks to the Lipschitz continuity of $\pi$, we infer that
\begin{eqnarray} \nonumber
- \frac{k\ell^2}{\nu} \int_{Q_t} [\pi(\phi_1) - \pi(\phi_2)](\phi_1 - \phi_2) 
& \leq & \frac{k\ell^2}{\nu} \int_{Q_t} |\pi(\phi_1) - \pi(\phi_2)| |\phi_1 - \phi_2| \\ \label{stimissimaamica3}
& \leq & \frac{k\ell^2 C_{\pi}}{\nu} \int_{Q_t} | \phi|^2. 
\end{eqnarray}
We  notice that the integral involving the gradients in \eqref{tutta1} is estimated from below in this way:
\begin{equation} \label{stimissimaamica4}
\int_{Q_t}  (|\nabla \Xi|^2 -\ell \nabla\phi\nabla\Xi + \ell^2 |\nabla\phi|^2) \geq 
\frac{1}{2} \int_{Q_t}  (|\nabla \Xi|^2 + \ell^2 |\nabla\phi|^2). 
\end{equation}
We also observe that 
\begin{equation} \label{stimissimaamica5}
-\frac{\gamma k\ell^3}{\nu} \int_{Q_t} |\phi|^2 \leq 0.
\end{equation}
Then, due to \eqref{stimissimaamica1}--\eqref{stimissimaamica5}, from \eqref{tutta1} we infer that
\begin{equation} \nonumber
\frac{1}{2} \int_{\Omega} |\Xi(t)|^2 + \frac{k\ell^2}{2\nu} \int_{\Omega} | \phi(t)|^2 + \frac{1}{2} \int_{Q_t}  (|\nabla \Xi|^2 + \ell^2 |\nabla\phi|^2)
\end{equation} 
\begin{equation} \nonumber
\leq \frac{1}{2} \| \Xi_0 \|^2_H
+ \frac{k\ell^2 C_{\pi}}{\nu} \int_{Q_t} |\phi|^2
+ \frac{k\ell^2}{2\nu} \| \phi_0 \|^2_H + \int_{Q_t} (f-k \Delta \Xi^{*})\Xi
+ \frac{\gamma k\ell^2}{\nu} \int_{Q_t} \Xi\phi 
+ \frac{\gamma k\ell^2}{\nu} \int_{Q_t}\Xi^{*}\phi .
\end{equation} 
By applying the inequality \eqref{dis1} to the last three terms of the right-hand side of the previous equation, we obtain that
\begin{equation} \nonumber
\frac{1}{2} \|\Xi(t)\|_H^2 + \frac{k\ell^2}{2\nu} \|\phi(t)\|_H^2 + \frac{1}{2} \int_{Q_t}  (|\nabla \Xi|^2 + \ell^2 |\nabla\phi|^2)
\end{equation} 
\begin{equation} \nonumber
\leq \frac{1}{2} \| \Xi_0 \|^2_H + \frac{k\ell^2}{2\nu} \| \phi_0 \|^2_H 
+ \frac{1}{8} \int_{Q_t} |\Xi|^2
+ 2 \int_{Q_t} |f-k \Delta \Xi^{*}|^2 + \frac{1}{8} \int_{Q_t} |\Xi|^2 
\end{equation}
\begin{equation} \label{precedentequi}
+ 2 \Big( \frac{\gamma k\ell^2}{\nu} \Big)^2 \int_{Q_t} |\phi|^2
+ \frac{1}{8} \int_{Q_t} |\Xi^{*}|^2 + 2\Big( \frac{\gamma k\ell^2}{\nu} \Big)^2 \int_{Q_t} |\phi|^2 + \frac{k\ell^2 C_{\pi}}{\nu} \int_{Q_t} |\phi|^2. 
\end{equation}
From \eqref{precedentequi} we infer that
\begin{equation} \nonumber
\frac{1}{2} \|\Xi(t)\|_H^2 + \frac{k\ell^2}{2\nu} \|\phi(t)\|_H^2 + \frac{1}{2} \int_{Q_t}  (|\nabla \Xi|^2 + \ell^2 |\nabla\phi|^2)
\end{equation} 
\begin{equation} \nonumber
\leq \frac{1}{2} \| \Xi_0 \|^2_H + \frac{k\ell^2}{2\nu} \| \phi_0 \|^2_H + 4 \| f \|^2_{L^2(Q)} + 4 k^2 T \| \Xi^{*} \|^2_W + \frac{1}{8}T \| \Xi^{*} \|^2_H
\end{equation} 
\begin{equation} \label{stimissimaamica7}
+ M \int^t_0 \Bigg( \frac{1}{2} \|\Xi (s)\|_H^2 +  \frac{k\ell^2}{2\nu}  \|\phi(s)\|_H^2 
+ \frac{1}{2} \int_{Q_s}  (|\nabla \Xi|^2 + \ell^2 |\nabla\phi|^2) \Bigg) \ ds,
\end{equation}
where
\begin{equation} \nonumber
M=\max{\Bigg(  \frac{ 4 \gamma^2 k \ell^2 + 2\nu C_{\pi}  }{\nu}  ; \frac{1}{2} \Bigg)}.
\end{equation}
From \eqref{stimissimaamica7}, by applying the Gronwall lemma, we conclude that
\begin{equation} \nonumber
\frac{1}{2} \|\Xi(t)\|_H^2 + \frac{k\ell^2}{2\nu} \|\phi(t)\|_H^2 
+ \frac{1}{2} \| \nabla \Xi \|^2_{L^{2}(0,t; H)} + \frac{\ell^2}{2} \| \nabla \phi \|^2_{L^{2}(0,t; H)}
\end{equation} 
\begin{equation} \label{stimissimaamica8}
\leq  C_1 \bigg[ 4 \|f \|^2_{L^2(Q)} + 4 k^2 T \| \Xi^{*} \|^2_W + \frac{1}{8} T \| \Xi^{*} \|^2_W  + C_0 \bigg( \| \Xi_0 \|^2_H + \| \phi_0 \|^2_H \bigg) \bigg], 
\end{equation}
where
\begin{equation} \nonumber
C_0 = \max{ \Bigg( \frac{1}{2} ; \ \frac{k\ell^2}{2\nu} \Bigg)}, \quad \quad \quad C_1= e^{TM}.
\end{equation}
From \eqref{stimissimaamica8}, we infer that 
\begin{equation} \nonumber 
C_3 \Big( \|\Xi(t) \|^2_H + \| \phi(t) \|_H^2 +  \| \nabla \Xi \|^2_{L^{2}(0,t; H)} + \| \nabla \phi \|^2_{L^{2}(0,t; H)} \Big)  
\end{equation}
\begin{equation} \nonumber
\leq \frac{1}{2} \|\Xi(t)\|_H^2 + \frac{k\ell^2}{2\nu} \|\phi(t)\|_H^2 
+ \frac{1}{2} \| \nabla \Xi \|^2_{L^{2}(0,t; H)} + \frac{\ell^2}{2} \| \nabla \phi \|^2_{L^{2}(0,t; H)}
\end{equation} 
\begin{equation} \nonumber
\leq C_2 \bigg( \|f \|^2_{L^2(Q)} + \| \Xi^{*} \|^2_W +  \| \Xi_0 \|^2_H + \| \phi_0 \|^2_H\bigg)
\end{equation}
\begin{equation} \label{eccoci22}
\leq C_2 \bigg( \|f \|_{L^2(Q)} + \| \Xi^{*} \|_W +  \| \Xi_0 \|_H + \| \phi_0 \|_H\bigg)^2,
\end{equation}
where
\begin{equation} \nonumber
C_2 = \max{ \Bigg( 4C_1 ; \ 4k^2T C_1 ; \ \frac{1}{8}TC_1 ; \ C_1 C_0  \Bigg)}, \quad \quad
C_3 = \min{\Bigg( \frac{1}{2} ; \ \frac{k\ell^2}{2\nu} ; \ \frac{\ell^2}{2} \Bigg)}.
\end{equation}
From \eqref{eccoci22} we obtain that
\begin{equation} \nonumber
\|\Xi(t) \|^2_H + \| \phi(t) \|^2_H +  \| \nabla \Xi \|^2_{L^{2}(0,t; H)} + \| \nabla \phi \|^2_{L^{2}(0,t; H)}
\end{equation}
\begin{equation} \label{quasifinita}
\leq C_4 \bigg( \|f \|_{L^2(Q)} + \| \Xi^{*} \|_W +  \| \Xi_0 \|_H + \| \phi_0 \|_H\bigg)^2,
\end{equation}
where $C_4 = \frac{C_2}{C_3}$. From \eqref{quasifinita} we conclude that
there exists a constant $C>0$ which depends only on $\Omega$, $T$ depends only on $\Omega$, $T$
and the parameters $\ell$, $\alpha$, $k$, $\nu$, $\gamma$ of the system, such that  
\begin{equation} \nonumber
\| \Xi_1 - \Xi_2 \|_{L^{\infty}(0,T; H) \cap L^{2}(0,T; V) } + \|  \phi_1 - \phi_2 \|_{L^{\infty}(0,T; H) \cap L^{2}(0,T; V) } 
\end{equation}
\begin{equation} \label{fineunicita1}
\leq C \big( \|f_1 -f_2 \|_{L^{2}(Q)} + \| \Xi^{*}_1 - \Xi^{*}_2 \|_W + \| \Xi_{0_1} - \Xi_{0_2} \|_H + \| \phi_{0_1} - \phi_{0_2} \|_H \big).
\end{equation}
To infer the uniqueness of the solution, we choose 
$f_1 = f_2$, $\Xi^*_1 = \Xi^*_2$, $\phi_{0_1} = \phi_{0_2}$, $\Xi_{0_1} = \Xi_{0_2}$.
Then, replacing the corresponding values in \eqref{fineunicita1}, we obtain that
\begin{equation} 
\| \Xi_1 - \Xi_2 \|_{L^{\infty}(0,T; H) \cap L^{2}(0,T; V)} + \|  \phi_1 - \phi_2 \|_{L^{\infty}(0,T; H) \cap L^{2}(0,T; V)} = 0.
\end{equation}
Hence $\Xi_1 = \Xi_2$ and $\phi_1 = \phi_2$.
Then the solution of problem $(P)$ (see \eqref{uni1}--\eqref{ultimafinale}) is unique.

\section*{Acknowledgments}
Partial support of the GNAMPA (Gruppo Nazionale per l'Analisi Matematica, 
la Pro\-babilit\`a e le loro Applicazioni) of INdAM is gratefully acknowledged.
The author is very grateful to the referees for the careful reading of the manuscript 
and for a number of useful suggestions.

\end{document}